\documentclass{gtmon_a}
\pdfoutput=1

\usepackage[dvips]{graphicx}
\usepackage[all]{xy}

\proceedingstitle{Exotic homology manifolds (Oberwolfach 2003)}
\conferencestart{29 June 2003}
\conferenceend{5 July 2003}
\conferencename{Workshop on Exotic Homology Manifolds}
\conferencelocation{Oberwolfach Mathematics Institute, Oberwolfach,
Germany}

\editor{Frank Quinn}
\givenname{Frank}
\surname{Quinn}

\editor{Andrew Ranicki}
\givenname{Andrew}
\surname{Ranicki}

\title{Controlled $L$--theory}

\author{Andrew Ranicki}
\givenname{Andrew}
\surname{Ranicki}
\address{School of Mathematics\\
University of Edinburgh\\\newline
Edinburgh\\
EH9 3JZ\\
UK}
\email{a.ranicki@ed.ac.uk}
\urladdr{}

\author{Masayuki Yamasaki}
\givenname{Masayuki}
\surname{Yamasaki}
\address{Department of Applied Science\\
Okayama University of Science\\\newline
Okayama\\
Okayama 700-0005\\
Japan}
\email{masayuki@mdas.ous.ac.jp}
\urladdr{}

\volumenumber{9}
\issuenumber{}
\publicationyear{2006}
\papernumber{7}
\startpage{105}
\endpage{153}

\doi{}
\MR{}
\Zbl{}

\arxivreference{math.GT/0402217}


\received{13 February 2004}
\revised{}
\accepted{13 February 2004}
\published{22 April 2006}
\publishedonline{22 April 2006}
\proposed{}
\seconded{}
\corresponding{}
\version{}

\makeatletter
\def\cnewtheorem#1[#2]#3{\newtheorem{#1}{#3}[section]
\expandafter\let\csname c@#1\endcsname\c@dummy}

\let\xysavmatrix\xymatrix
\def\xymatrix{\disablesubscriptcorrection\xysavmatrix}
\AtBeginDocument{\let\wbar\overline\let\wtilde\widetilde\let\hat\widehat}

\UseComputerModernTips

\cnewtheorem{theorem}[dummy]{Theorem}
\cnewtheorem{lemma}[dummy]{Lemma}
\cnewtheorem{corollary}[dummy]{Corollary}
\cnewtheorem{proposition}[dummy]{Proposition}

\cnewtheorem{conjecture}[dummy]{Conjecture}
\theoremstyle{definition}               
\cnewtheorem{definition}[dummy]{Definition}
\cnewtheorem{question}[dummy]{Question}
\cnewtheorem{example}[dummy]{Example}

\newtheorem*{remarkun}{Remark}
\newtheorem*{remarksun}{Remarks}

\makeatother  

\newcommand{\ZZ}{\mathbb Z}
\newcommand{\RR}{\mathbb R}
\newcommand{\QQ}{\mathbb Q}
\newcommand{\LL}{\mathbb L}
\renewcommand{\C}{\mathcal C}
\renewcommand{\D}{\mathcal D}
\newcommand{\F}{\mathcal F}
\def\exc{{\rm {exc}}}
\def\ra#1{\hbox to #1pc{\rightarrowfill}}
\def\la#1{\hbox to #1pc{\leftarrowfill}}
\def\sa{\ra{1.5}}
\def\lga{\ra{2.5}}
\def\d{{\partial}}
\long\def\showlabel#1{}

\makeop{Wh}

\begin{document}

\begin{asciiabstract}
We develop an epsilon-controlled algebraic L-theory, extending our
earlier work on epsilon-controlled algebraic K-theory.  The controlled
L-theory is very close to being a generalized homology theory; we study
analogues of the homology exact sequence of a pair, excision properties,
and the Mayer--Vietoris exact sequence.  As an application we give a
controlled L-theory proof of the classic theorem of Novikov on the
topological invariance of the rational Pontrjagin classes.
\end{asciiabstract}

\begin{abstract}
We develop an epsilon-controlled algebraic $L$--theory, extending our
earlier work on epsilon-controlled algebraic $K$--theory.  The controlled
$L$--theory is very close to being a generalized homology theory; we study
analogues of the homology exact sequence of a pair, excision properties,
and the Mayer--Vietoris exact sequence.  As an application we give a
controlled $L$--theory proof of the classic theorem of Novikov on the
topological invariance of the rational Pontrjagin classes.
\end{abstract}

\maketitle

\section{Introduction}\label{intro}

\showlabel{intro}

The purpose of this article is to develop a controlled algebraic $L$--theory,
of the type first proposed by Quinn \cite{Q3} in connection with the
resolution of homology manifolds by topological manifolds.

We define and study the epsilon-controlled $L$--groups
$L_n^{\delta,\epsilon}(X;p_X,R)$, extending to $L$--theory the controlled
$K$--theory of Ranicki and Yamasaki \cite{RY-K}.  When the control map
$p_X$ is a fibration and $X$ is a compact ANR, these groups are stable in
the sense that $L_n^{\delta,\epsilon}(X;p_X,R)$ depends only on $p_X$ and
$R$ and not on $\delta$ or $\epsilon$ as long as $\delta$ is sufficiently
small and $\epsilon \ll \delta$ (see Pedersen--Yamasaki \cite{P-Y}).

These are the candidates of the controlled surgery obstruction groups;
in fact, such a controlled surgery theory has been established when the
control map $p_X$ is $UV^1$ (see Pedersen--Quinn--Ranicki \cite{PQR}
and Ferry \cite{Ferry}).

Although epsilon controlled $L$--groups do not produce a homology
theory in general, they have the features of a generalized
homology modulo controlled $K$--theory problems. In this article we
study the controlled $L$--theory analogues of the homology exact
sequence of a pair (\fullref{ThESP}), excision properties
(\fullref{Exc}), and the Mayer--Vietoris sequence (\fullref{ThMVS1}).

In certain cases when there are no controlled $K$--theoretic
difficulties, we can actually show that controlled $L$--groups are
generalized homology groups. This is discussed in \fullref{SpC}.

In the last two sections, we study locally-finite analogues and as an
application give a controlled $L$--theory proof of the classic theorem
of Novikov \cite{Novikov} on the topological invariance of the rational
Pontrjagin classes.

\section{Epsilon-controlled quadratic structures}\label{q-str}

\showlabel{q-str}
In this section we study several operations concerning quadratic Poincar\'e
complexes with geometric control.  These will be used to define epsilon controlled
$L$--groups in the next section.

In \cite{RY-K} we discussed various aspects of geometric modules and
morphisms and geometric control on them, and studied $K$--theoretic
properties of geometric (=free) and projective module chain complexes
with geometric control.  There we considered only $\ZZ$--coefficient
geometric modules, but the material in Sections \ref{intro}--\ref{MV}
remains valid if we use any ring $R$ with unity as the coefficient.
To incorporate the coefficient ring into the notation, the group
$\wtilde K_0(X,p_X,n,\epsilon)$ defined using the coefficient ring $R$
will be denoted $\wtilde K_0^{n,\epsilon}(X; p_X,R)$ in this article.

To deal with $L$--theory, we need to use duals. Fix the control map
$p_X\co M\to X$ from a space $M$ to a metric space $X$ and let $R$ be
a ring with involution (see Ranicki \cite{AdditiveL}). The \emph{dual} $G^*$ of a
geometric $R$--module $G$ is $G$ itself. Recall that a geometric
morphism is a linear combination of paths in $M$ with coefficient
in $R$. The \emph{dual} $f^*$ of a geometric morphism
$f=\sum_\lambda a_\lambda \rho_\lambda$ is defined by
$f^*=\sum_\lambda \wbar a_\lambda \wbar\rho_\lambda$, where $\wbar
a_\lambda\in R$ is the image of $a$ by the involution of $R$ and
$\wbar\rho_\lambda$ is the path obtained from $\rho_\lambda$ by
reversing the orientation. Note that if $f$ has radius $\epsilon$
then so does its dual $f^*$ and that $f\sim _{\epsilon} g$ implies
$f^* \sim _{\epsilon}g^*$, by our convention. For a geometric
$R$--module chain complex $C$, its $n${\sl-dual} $C^{n-*}$ is
defined using the formula in Ranicki \cite{Exact}.

For a subset $S$ of a metric space $X$, $S^\epsilon$ will denote the closed
$\epsilon$ neighborhood of $S$ in $X$ when $\epsilon\ge0$.  When $\epsilon<0$,
$S^{\epsilon}$ will denote the set $X-(X-S)^{-\epsilon}$.

Let $C$ be a free $R$--module chain complex on $p_X\co M\to X$.
An $n$--\emph{dimensional $\epsilon$ quadratic structure} $\psi$ \emph{on} $C$
is a collection
$\{\psi_s | s\ge 0\}$ of geometric morphisms
\[
\psi_s \co  C^{n-r-s} = (C_{n-r-s})^* \to C_r \qquad \text{($r\in \ZZ$)}
\]
of radius $\epsilon$ such that
\begin{equation}
\label{eqn:star}
d\psi_s + (-)^r\psi_sd^* +(-)^{n-s-1}(\psi_{s+1}+(-)^{s+1}T\psi_{s+1})
\sim_{3\epsilon} 0\co C^{n-r-s-1} \to C_r, \tag{$*$}
\end{equation}
for $s\ge0$.
An $n$--dimensional free $\epsilon$ chain complex $C$ on $p_X$ equipped with an
$n$--dimensional $\epsilon$ quadratic structure is called an
$n$--\emph{dimensional $\epsilon$ quadratic $R$--module complex on} $p_X$.

Let $f\co C\to D$ be a chain map between free
chain complexes on $p_X$.
An $(n{+}1)$--\emph{dimensional $\epsilon$ quadratic structure $(\delta\psi,\psi)$
on} $f$ is a collection $\{\delta\psi_s,\psi_s |s\ge 0\}$ of geometric
morphisms $\delta\psi_s\co D^{n+1-r-s} \to D_r$, $\psi_s\co C^{n-r-s} \to C_r$
($r \in \ZZ$) of radius $\epsilon$ such that the following holds in
addition to \eqref{eqn:star}:
\begin{multline*}
d(\delta\psi_s)+(-)^r(\delta\psi_s)d^*+(-)^{n-s}\bigl(\delta\psi_{s+1}
+(-)^{s+1}T\delta\psi_{s+1}\bigr) +(-)^nf\psi_s f^* \\
\sim_{3\epsilon} 0 \co D^{n-r-s} \to D_r,
\end{multline*}
for $s\ge0$.
An $\epsilon$ chain map $f\co C\to D$ between an $n$--dimensional free $\epsilon$
chain complex $C$ on $p_X$ and an $(n{+}1)$--dimensional free $\epsilon$ chain
complex $D$ on $p_X$ equipped with an $(n{+}1)$--dimensional $\epsilon$ quadratic
structure is called an $(n{+}1)$--\emph{dimensional $\epsilon$ quadratic
$R$--module pair on} $p_X$.
Obviously its \emph{boundary} $(C,\psi)$ is an $n$--dimensional $\epsilon$
quadratic $R$--module complex on $p_X$.
We will suppress references to the coefficient ring $R$ unless we need to emphasize the coefficient ring.

An $\epsilon$ \emph{cobordism} of $n$--dimensional $\epsilon$ quadratic
structures $\psi$ on $C$ and $\psi'$ on $C'$ is an $(n{+}1)$--dimensional
$\epsilon$ quadratic structure $(\delta\psi,\psi\oplus -\psi')$ on some
chain map $C\oplus C'\to D$.
An $\epsilon$ \emph{cobordism} of $n$--dimensional $\epsilon$ quadratic complexes
$(C,\psi)$, $(C',\psi')$ on $p_X$ is an $(n{+}1)$--dimensional $\epsilon$
quadratic pair on $p_X$
\[
\bigl(\begin{pmatrix}f & f'\end{pmatrix} \co  C\oplus C' \to
D,\quad(\delta\psi,\psi\oplus -\psi')\bigr)
\]
with boundary $(C\oplus C',\psi\oplus -\psi')$. The union of
adjoining cobordisms is defined using the formula in Ranicki \cite{Exact}.
The union of adjoining $\epsilon$ cobordisms is a $2\epsilon$
cobordism.

$\Sigma C$ and $\Omega C$ will denote the suspension and the desuspension of
$C$ respectively, and $\C(f)$ will denote the algebraic mapping cone
of a chain map $f$.

\begin{definition}
Let $W$ be a subset of $X$.
An $n$--dimensional $\epsilon$ quadratic structure $\psi$ on $C$ is
$\epsilon$ \emph{Poincar\'e (over $W$)} if
the algebraic mapping cone of the duality $3\epsilon$ chain map
\[
\D_{\psi}=(1+T)\psi_0\co C^{n-*} \lga C
\]
is $4\epsilon$ contractible (over $W$).
A quadratic complex $(C,\psi)$ is
$\epsilon$ \emph{Poincar\'e (over $W$)} if  $\psi$ is $\epsilon$ Poincar\'e
(over $W$).
Similarly, an $(n{+}1)$--dimensional $\epsilon$ quadratic structure
$(\delta\psi,\psi)$ on $f\co C\to D$ is $\epsilon$ \emph{Poincar\'e (over $W$)}
if the algebraic mapping cone of the duality $4\epsilon$ chain map
\[
\D_{(\delta\psi,\psi)}=
\begin{pmatrix}(1+T)\delta\psi_0 \cr
(-)^{n+1-r}(1+T)\psi_0f^*\cr\end{pmatrix}\co 
D^{n+1-r} \to \C(f)_r = D_r\oplus C_{r-1}
\]
is $4\epsilon$ contractible (over $W$) (or equivalently the algebraic mapping
cone of the $4\epsilon$ chain map
\[
\wbar\D_{(\delta\psi,\psi)}=((1+T)\delta\psi_0\quad f(1+T)\psi_0)\co \C(f)^{n+1-*} \lga D
\]
is $4\epsilon$ contractible (over $W$)) and $\psi$ is $\epsilon$
Poincar\'e (over $W$).
A quadratic pair $(f,(\delta\psi,\psi))$ is \emph{$\epsilon$ Poincar\'e
(over $W$)} if $(\delta\psi,\psi)$ is $\epsilon$ Poincar\'e (over $W$).
We will also use the notation $\D_{\delta\psi}=(1+T)\delta\psi_0$, although
it does not define a chain map from $D^{n+1-*}$ to $D$ in general.
\end{definition}

This definition is slightly different from the one given in Yamasaki \cite{MY87}
(especially when $W$ is a proper subset of $X$).
There a quadratic complex/pair was defined to be
$\epsilon$ Poincar\'e over $W$ if the duality map is an $\epsilon$ chain
equivalence over $W$.
If $\C(\D_{\psi})$ (resp. $\C(\D_{(\delta\psi,\psi)})$) is $4\epsilon$
contractible (over $W$), then
$ \D_{\psi}$ (resp. $\D_{(\delta\psi,\psi)}$ ) is only a ``weak''
$8\epsilon$ chain equivalence over $W$.

\begin{definition}
\label{weak-eq}
\showlabel{weak-eq}
A chain map $f\co C \to D$ is a \emph{weak $\epsilon$ chain equivalence
over} $W$ if
\begin{enumerate}
\item $f$ is an $\epsilon$ chain map;
\item there exists a family $g=\{g_r\co D_r\to C_r\}$ of geometric morphisms
of radius $\epsilon$ such that
\begin{itemize}
\item $dg_r$ and $g_rd$ have radius $\epsilon$, and
\item $dg_r ~\sim_\epsilon ~g_{r-1}d$ over $W$
\end{itemize}
for all $r$;
\item there exist two families $h=\{h_r\co C_r\to C_{r+1}\}$ and
$k=\{k_r\co D_r\to D_{r+1}\}$ of $\epsilon$ morphisms such that
\begin{itemize}
\item $dh_r + h_{r-1}d ~\sim_{2\epsilon}~ 1-g_rf_r$ over $W$, and
\item $dk_r + k_{r-1}d ~\sim_{2\epsilon}~ 1-f_rg_r$ over $W$.
\end{itemize}
for all $r$.
\end{enumerate}
In other words a weak chain equivalence satisfies all the
properties of a chain equivalence except that its inverse may not
be a chain map outside of $W$.
\end{definition}

Weak chain equivalences behave quite similarly to chain equivalences.
For example, 2.3(3) and 2.4 of Ranicki \cite{RY-K} can be easily
generalized as follows:

\begin{proposition}
\label{comp}
\showlabel{comp}
If $f\co C\to D$ is a weak $\delta$ chain equivalence over $V$ and
$f'\co D\to E$ is a weak $\epsilon$ chain equivalence over $W$,
then $f'f$ is a weak $\delta+\epsilon$ chain equivalence over
$V^{-\delta-\epsilon} \cap W^{-\delta}$.  If we further assume that
$f$ is a $\delta$ chain equivalence, then $f'f$ is a weak $\delta+\epsilon$
chain equivalence over $V^{-\epsilon} \cap W^{-\delta}$.
\end{proposition}

\begin{proposition}
\label{cone}
\showlabel{cone}
Let $f\co C\to D$ be an $\epsilon$ chain map.
If the algebraic mapping cone $\C(f)$ is $\epsilon$ contractible over $W$,
then $f$ is a weak $2\epsilon$ chain equivalence over $W$.  If $f$ is a
weak $\epsilon$ chain equivalence over $W$, then $\C(f)$ is $3\epsilon$
contractible over $W^{-2\epsilon}$.
\end{proposition}

We have employed the definition of Poincar\'e complexes/pairs using local
contractibility of the algebraic mapping cone of the duality map,
because algebraic mapping cones are easier to handle than chain equivalences.
For example, consider
a triad $\Gamma$:
\[
\xymatrix@R=5pt{ C \ar[dd]_{g} \ar[r]^{f} \ar @{~>} [ddr]^{h} &  D \ar[dd]^{g'}\\
& & {dh + hd\sim f'g - g'f}\\
C' \ar[r]_{f'} & D'
}
\]
and assume
\begin{enumerate}
\item $f$ (resp. $f'$) is a $\delta$ (resp. $\delta'$) chain map,
\item $g$ (resp. $g'$) is an $\epsilon$ (resp. $\epsilon'$) chain map,
\item $h\co g'f\simeq f'g $ is a $\gamma$ chain homotopy.
\end{enumerate}
Then there are induced a $\max\{\delta,\delta',2\gamma\}$ chain map
\[
F=\begin{pmatrix} f' & (-)^rh \\ 0 & -f \cr\end{pmatrix} \co 
\C(-g)_r=C'_r\oplus C_{r-1} \to \C(g')_r=D'_r\oplus D_{r-1}
\]
and a $\max\{\epsilon,\epsilon',2\gamma\}$ chain map
\[
G=\begin{pmatrix}g' & (-)^rh \\ 0 & g \end{pmatrix} \co 
\C(f)_r=D_r\oplus C_{r-1} \to \C(f')_r=D'_r\oplus C'_{r-1}.
\]
It is easily seen that $\C (F) = \C (G)$.
%

\begin{proposition}
\label{cone-g}
\showlabel{cone-g}
If $\C(g\co C\to C')$ is $\epsilon$ contractible
over $W$, then $\C(-g)$ is $\epsilon$ contractible over $W$.
\end{proposition}

\begin{proof}
Suppose
\[
\begin{pmatrix}a & b \\ c & d \end{pmatrix}\co\C(g)_r=C_r\oplus C'_{r-1} ~\sa~
\C(g)_{r+1}=C_{r+1}\oplus C'_r
\]
is an $\epsilon$ chain contraction over $W$ of $\C(g)$, then
\[
\begin{pmatrix}a & -b \\ -c & d\end{pmatrix}
\]
is an $\epsilon$ chain contraction over $W$ of $\C(-g)$.
\end{proof}

\begin{proposition}
\label{cont-triad}
\showlabel{cont-triad}
Let $\Gamma$ be as above, and further assume that
$\C(g)$ is $\epsilon$ contractible over $W$ and $\C(g')$ is $\epsilon'$
contractible over $W$, then $\C(G)$ is
$3\max\{\epsilon,\epsilon',\delta,\delta',2\gamma\}$ contractible
over $W^{-2\max\{\epsilon,\epsilon',\delta,\delta',2\gamma\}}$.
\end{proposition}

\begin{proof}
By \fullref{cone-g}, $\C(-g)$ is $\epsilon$ contractible over $W$.
Therefore $F\co \C(-g) \to \C(g')$ is a
$\max\{\epsilon, \epsilon', \delta, \delta', 2\gamma\}$
chain equivalence over $W$, and the proposition is proved by
applying \fullref{cone} to $F$.
\end{proof}

\begin{corollary}
\label{P-P}
\showlabel{P-P}
Let $C$ and $D$ be free $\epsilon$ chain complexes,
and let $(\delta\psi,\psi)$ be an $\epsilon$ quadratic structure on
an $\epsilon$ chain map $f\co C\to D$.  If $\C(\D_{(\delta\psi,\psi)})$ is
$4\epsilon$ contractible over $W$, then $\C(\D_\psi)$ is $100\epsilon$
contractible over $W^{-100\epsilon}$.
\end{corollary}

\begin{proof}
Consider the triad $\Gamma$:
\[
\xymatrix@R=5pt{\Omega\C(f)^{n+1-*}\ar[dd]_{\D_{(\delta\psi,\psi)}} \ar[rr]^{\alpha=(1~~0)}
\ar @{~>} [ddrr]^{h} & & \Omega\D^{(n+1-*)} \ar[dd]^{\wbar\D_{(\delta\psi,\psi)}}\\
& & &  {h=\begin{pmatrix}0 & 0\\ 0 & (-)^{r+1}(1+T)\psi_0\end{pmatrix}}\\
\Omega\D \ar[rr]_{\beta={}^t(1~~0)} & & \Omega\C(f)}
\]

\noindent
and consider the chain map $G\co \C(\alpha) \to \C(\beta)$ induced from
$\Gamma$ as above.  Then $\C(G)$ is $12\epsilon$ contractible
over $W^{-8\epsilon}$ by the previous proposition. Therefore
$G$ is a weak $24\epsilon$ chain equivalence over $W^{-8\epsilon}$.
$(1+T)\psi_0$ is equal to the following composite of $G$ with two $\epsilon$
chain equivalences:
\[
C^{n-*} \xrightarrow[\simeq_{\epsilon}]{{}^t(0\,\,0\,\,1)}
\C(\alpha) \xrightarrow{~~~G~~~}
\C(\beta) \xrightarrow[\simeq_{\epsilon}]{(0\,\,1\,\,0)}  C,
\]
and the claim follows from  \fullref{comp}.
\end{proof}

Next we describe various constructions on quadratic complexes
with some size estimates.  Firstly
a direct calculation shows the following.  (See the non-controlled case
in Ranicki \cite{Exact}.)

\begin{proposition}
\label{union}
\showlabel{union}
If adjoining $\epsilon$ cobordisms $c$ and $c'$
are $\epsilon$ Poincar\'e over $W$, then $c\cup c'$ is $100\epsilon$
Poincar\'e over $W^{-100\epsilon}$.
\end{proposition}

The following proposition gives us a method to construct quadratic
structures and cobordisms.

\begin{proposition}
\label{i-str}
\showlabel{i-str}
Suppose
$g\co C\to C'$ is a $\delta$ chain map of $\delta$ chain complexes and
$\psi$ is an $n$--dimensional $\epsilon$ quadratic structure on $C$.
\par\noindent{\rm (1)}\qua
$g_\%\psi = \{ (g_\%\psi)_s = g\psi_sg^*\}$ is a
$2\delta+\epsilon$ quadratic structure on $C'$, and
$(0,\psi\oplus -g_\%\psi)$
is a $2\delta+\epsilon$ quadratic structure on the $\delta$ chain map
$(g\quad 1)\co C\oplus C' \sa C'$.
\par\noindent{\rm (2)}\qua
If $\psi$ is $\epsilon$ Poincar\'e over $W$ and $g$ is
a weak $\delta$ chain equivalence over $W$,
then $g_\%\psi$ and $(0,\psi\oplus -g_\%\psi)$ are $2\delta+6\epsilon$ Poincar\'e over
$W^{-(6\delta+24\epsilon)}$.
\par\noindent{\rm (3)}\qua
If $\psi$ is $\epsilon$ Poincar\'e, $g$ is a $\delta$ chain
equivalence, $\Psi'=(\delta\psi,\psi'=g_\%\psi)$ is an $\epsilon'$
Poincar\'e $\delta'$ quadratic structure on a $\delta'$ chain map
$f'\co C'\to D$, and $D$ is a $\gamma$ chain complex,
then $\Psi=(\delta\psi,\psi)$ is an $\epsilon'+3\max\{9\delta,
6\delta',4\epsilon,3\gamma\}$ Poincar\'e $\max\{\delta,\epsilon\}$
quadratic structure on the $\delta'+\epsilon$ chain map $f=f'\circ g\co C\to D$.
\end{proposition}

\begin{proof}
(1)\qua This can be checked easily.

(2)\qua This holds because the duality maps for $(C',g_\%\psi)$ and $c$
split as follows:
\begin{gather*}
(C')^{n-*} \xrightarrow{~~g^*~~}  ~~C^{n-*}~~ \xrightarrow{(1+T)\psi_0}
C \xrightarrow {~~g~~}  C'\\
\C((g\,\,\,\,1))^{n+1-*}~~ \xrightarrow[\simeq_{\delta}] {(0\,\,\,\,1\,\,-g^*)}
~~C^{n-*} \xrightarrow{(1+T)\psi_0} C \xrightarrow{~~g~~} C'
\end{gather*}
(3)\qua We study the duality map for $\Psi'$.
Since $g$ is a $\delta$ chain equivalence and $\C(1\co D\to D)$ is
$\gamma$ contractible, the algebraic mapping cone of the $\max\{\gamma,
\delta\}$ chain map
\[
\tilde g = \begin{pmatrix}1 & 0 \\ 0 & g\end{pmatrix}\co \C(f)_r=D_r\oplus C_{r-1}\lga \C(f')=
D_r\oplus C'_{r-1}
\]
is $3\max\{3\delta, \delta+\delta',\gamma\}$ contractible, and
so is $\C({\tilde g}^*\co  \C(f')^{n+1-*}\to \C(f)^{n+1-*})$.
Therefore, the chain map $\C(\D_{\Psi'})\to \C(\D_\Psi)$ defined by
\[
\begin{pmatrix}1 & 0\\ 0 & {\tilde g}^*\end{pmatrix}\co \C(\D_{\Psi'})=D_r\oplus\C(f')^{n+2-r}
\lga\C(\D_{\Psi})=D_r\oplus\C(f)^{n+2-r}
\]
is a $6\max\{9\delta,6\delta',4\epsilon,3\gamma\}$ chain equivalence.
The claim now follows from the next lemma.
\end{proof}

\begin{lemma}
\label{l-cont}
\showlabel{l-cont}
If a chain complex $A$ is $\epsilon$ chain equivalent to a chain complex $B$
which is $\delta$ contractible over $X-Y$, then $A$ is
$(2\epsilon+\delta)$ contractible over $X-Y^\epsilon$.
\end{lemma}

\begin{proof}
Let $f\co A\to B$ be an $\epsilon$ chain equivalence,
$g$ an $\epsilon$ chain homotopy inverse, $h\co gf\simeq_\epsilon 1$ an $\epsilon$
chain homotopy, and $\Gamma$ a $\delta$ chain contraction of $B$ over $X-Y$.
Then $g\Gamma f+h$ gives a$(2\epsilon+\delta)$ chain contraction of $A$
over $X-Y^\epsilon$.
\end{proof}

\begin{remarksun}
(1)\qua An $\epsilon$ chain equivalence $f\co C\sa C'$ such that
$f_\%\psi=\psi'$ will be called an $\epsilon$ \emph{homotopy equivalence}
from $(C,\psi)$ to $(C',\psi')$.  By \fullref{i-str}, a homotopy equivalence between
quadratic Poincar\'e complexes induces a Poincar\'e cobordism between them.

(2)\qua The estimates given in  \fullref{i-str} and
\fullref{l-cont} are, of course, not acute in general.
For example, consider an $\epsilon$ quadratic complex $(C,\psi)$
which is $\epsilon$ Poincar\'e over $W$.  Then a direct calculation shows that the cobordism
between $(C,\psi)$ and itself induced by the identity map
of $C$ is an $\epsilon$ quadratic pair and is $\epsilon$ Poincar\'e
over $W$.  This cobordism will be called the \emph{trivial cobordism}
from $(C,\psi)$ to itself.
\end{remarksun}

\section{Epsilon-controlled $L$--groups}\label{L-groups}

\showlabel{L-groups} In this section we review the boundary
construction of the first-named author and then introduce
epsilon-controlled $L$--groups
\[
L_n^{\delta,\epsilon}(X;p_X,R)  \quad \text{and} \quad L_n^{\delta,\epsilon}(X,Y;p_X,R)
\]
for $p_X\co M\to X$, $Y\subset X$, $n\ge 0$, $\delta\ge\epsilon\ge0$, and a ring $R$ with involution.
These are defined using geometric $R$--module chain complexes with quadratic
Poincar\'e structures discussed in the previous section.

Let $(C,\psi)$ be an $n$--dimension\-al $\epsilon$ quadratic $R$--module complex on $p_X$,
where $n\ge 1$.
Define a (possibly non-positive) $2\epsilon$ chain complex $\d C$
by $\Omega\C(\D_{\psi})$.
Then an $(n{-}1)$--dimensional $2\epsilon$ Poincar\'e $2\epsilon$ quadratic
structure $\d\psi$ on $\d C$ is defined by:
\begin{align*}
\d\psi_0&=\begin{pmatrix}0 & 0 \\ 1 & 0\end{pmatrix}\co \d C^{n-r-1}=C^{n-r}\oplus C_{r+1}
\sa \d C_r=C_{r+1}\oplus C^{n-r}\\
\d\psi_s &=\begin{pmatrix}(-)^{n-r-s-1}T\psi_{s-1} & 0\\
0&0\end{pmatrix}\co \\
&\quad \d C^{n-r-s-1}=C^{n-r-s}\oplus C_{r+s+1} \sa\d C_r=C_{r+1}\oplus C^{n-r}
\qquad (s\ge 1).\cr
\end{align*}
This is the \emph{boundary construction} of Ranicki \cite{Exact}.
The structure $\Psi_1=(0,\d\psi)$ is an $n$--dimensional $2\epsilon$ Poincar\'e
$2\epsilon$ quadratic structure on the $\epsilon$ chain map
\[
i_C={\rm projection}\co \d C \sa C^{n-*}
\]
of $2\epsilon$ chain complexes.
This is called the \emph{algebraic Poincar\'e thickening} (see Ranicki
\cite{Exact}).

\begin{example}
\label{b-ex}
\showlabel{b-ex}
Consider an $n$--dimensional $\epsilon$ chain complex $F$, and
give $\Sigma F$ the trivial $(n{+}1)$--dimensional quadratic structure
$\theta_s = 0$ ($s\ge 0$).  Its algebraic Poincar\'e thickening
\[
(i_{\Sigma F} \co  \d\Sigma F = F\oplus F^{n-*} \sa (\Sigma F)^{n+1-*}=
F^{n-*}, \quad (0,\d\theta))
\]
is an $(n{+}1)$--dimensional $\epsilon$ Poincar\'e $\epsilon$ null-cobordism
of $(\d\Sigma F, \d\theta)$.
\end{example}

There is an inverse operation up to homotopy equivalence.
Given an $n$--dimension\-al $\epsilon$ Poincar\'e $\epsilon$
quadratic pair $c=(f\co C\to D,(\delta\psi,\psi))$, take the union
$(\tilde C,\tilde\psi)$ of $c$ with the $\epsilon$ quadratic pair
$(C\to 0,(0,-\psi))$.
$\tilde C$ is equal to $\C(f)$.
$(\tilde C,\tilde \psi)$ is an $n$--dimensional $2\epsilon$ quadratic complex
and is called the \emph{algebraic Thom complex} of $c$.
The algebraic Poincar\'e thickening of $(\tilde C,\tilde\psi)$ is
``homotopy equivalent'' to the original pair $c$ (as pairs).
Since we will not use this full statement, we do not define homotopy
equivalences of pairs here and only mention that the chain map
\[
g=\begin{pmatrix}0&1&0&-\psi_0\end{pmatrix}\co \d\tilde C_r=D_{r+1}\oplus C_r\oplus D^{n-r}
\oplus C^{n-r-1} \sa C_r
\]
gives an $11\epsilon$ chain equivalence such that
$g_\%(\d\tilde\psi)=\psi$.
If we start with an $n$--dimensional $\epsilon$ quadratic complex $(C,\psi)$
on $p_X$, then the algebraic Thom complex of the algebraic Poincar\'e
thickening $(i_C\co \d C\sa C^{n-*}, (0,\d\psi))$ of $(C,\psi)$ is
$3\epsilon$ homotopy equivalent to $(C,\psi)$;
$3\epsilon$ homotopy equivalences are given by
\begin{gather*}
f=\begin{pmatrix}-\D_\psi&1&0\end{pmatrix}\co \C(i_C)_r=C^{n-r}\oplus C_r\oplus C^{n-r+1}\sa C_r,\\
f_\%(\Psi_1\cup_{\d\psi}-\Psi_2)=\psi,\\
f'={}^t\begin{pmatrix}0&1&0\end{pmatrix}\co C_r\sa\C(i_C)_r=C^{n-r}\oplus C_r\oplus C^{n-r+1},\\
f'_\%\psi=\Psi_1\cup_{\d\psi}-\Psi_2,
\end{gather*}
where $\Psi_2=(0,\d\psi)$ is the $n$--dimensional $\epsilon$ quadratic
structure on the trivial chain map $0\co \d C\sa 0$.

The boundary construction described above generalizes to quadratic pairs.
For an $(n{+}1)$--dimensional $\epsilon$ quadratic pair $(f\co C\to D,(\delta\psi,\psi))$ on $p_X$,
define a (possibly non-positive) $2\epsilon$ chain complex
 $\d D$ by $\Omega\C(\D_{(\delta\psi,\psi)})$ and define
an $n$--dimensional $3\epsilon$ Poincar\'e $2\epsilon$
quadratic structure $\Psi_3=(\d\delta\psi,\d\psi)$ on the $2\epsilon$ chain map
of $2\epsilon$ chain complexes
\[
\d f=\begin{pmatrix}f&0\\0&0\\0&1\end{pmatrix}\co \d C_r=C_{r+1}\oplus C^{n-r} \lga
\d D_r=D_{r+1}\oplus D^{n-r+1}\oplus C^{n-r}
\]
by
\begin{align*}
\d\delta\psi_0 &=\begin{pmatrix}0&0&0\\-1&0&0\\0&0&0\end{pmatrix}\co 
\d D^{n-r}=D^{n-r+1}\oplus D_{r+1}\oplus C_r \\
&\qquad\qquad\qquad\qquad\lga \d D_r=D_{r+1}\oplus D^{n-r+1}\oplus C^{n-r}\\
\d\delta\psi_s &=\begin{pmatrix}(-)^{n-r-s-1}T\delta\psi_{s-1}&0&0\\0&0&0\\0&0&0\end{pmatrix}
\co \d D^{n-r-s}=D^{n-r-s+1}\oplus D_{r+s+1}\oplus C_{r+s}\\
&\qquad\qquad\qquad\qquad\lga \d D_r=D_{r+1}\oplus D^{n-r+1}\oplus C^{n-r}\qquad (s\ge 1).
\end{align*}
$\d\psi$ is the same as above.
Then $(0,\Psi_1\cup_{\d\psi}(-\Psi_3))$ gives an $(n{+}1)$--dimensional
$300\epsilon$ Poincar\'e $4\epsilon$ quadratic structure on the $\epsilon$
chain map
\begin{align*}
\begin{pmatrix}i&0&i_D\end{pmatrix} &\co (C^{n-*}\cup_{\d C}\d D)_r=C^{n-r}\oplus \d C_{r-1}\oplus\d D_r\\
&\quad \lga (\C(f)^{n+1-*})_r=\C(f)^{n+1-r}
\end{align*}
of $2\epsilon$ chain complexes,
where $i\co C^{n-*}\to \C(f)^{n+1-*}$ is the inclusion map and
$i_D\co \d D\to\C(f)^{n+1-*}$ is the projection map.

If $(C,\psi)$ (resp. $(f\co C\to D,(\delta\psi,\psi))$) is $\epsilon$
Poincar\'e, then $\d C$ is (resp. $\d C$ and $\d D$ are)
$4\epsilon$ contractible, and hence chain homotopic to a positive
chain complex (resp. positive chain complexes). But in general $\d
C$ (and $\d D$) may not be chain homotopic to a positive chain
complex. This leads us to the following definition. The
non-controlled version is described in Ranicki \cite{Exact}.

\begin{definition} (1)\qua A positive geometric chain complex $C$ ($C_i=0$
for $i<0$) is $\epsilon$ \emph{connected} if there exists a $4\epsilon$ morphism
$h\co C_0 \to C_{1}$  such that $dh \sim_{8\epsilon}1_{C_0}$.

(2)\qua A chain map $f\co C\to D$ of positive chain complexes is
$\epsilon$ \emph{connected} if $\C(f)$ is $\epsilon$ connected.

(3)\qua A quadratic complex $(C,\psi)$ is $\epsilon$ \emph{connected} if
$\D_{\psi}$ is $\epsilon$ connected.

(4)\qua A quadratic pair $(f\co C\to D,(\delta\psi,\psi))$ is $\epsilon$
\emph{connected} if $\D_{\psi}$ and $\D_{(\delta\psi,\psi)}$ are $\epsilon$
connected.
\end{definition}

\begin{lemma}
\label{conn}
\showlabel{conn}
{\rm (1)}\qua The composition of a $\delta$ connected chain map and an $\epsilon$ connected chain map
is $\delta+\epsilon$ connected.

{\rm (2)}\qua Quadratic complexes and pairs that are $\epsilon$ Poincar\'e are $\epsilon$ connected.

{\rm (3)}\qua If $\psi$ is an $\epsilon$ connected quadratic structure
on $C$ and $g\co C\to C'$ is a $\delta$ connected chain map,
then $\D_{(0,\psi\oplus -g_\%\psi)}$ is $\epsilon+2\delta$ connected.
\end{lemma}

\begin{proof} (1) is similar to \fullref{comp}.  (2) is immediate from definition.
(3) is similar to \fullref{i-str} (2).
\end{proof}

\begin{remarkun}
In general the $\epsilon$ connectivity of $g$ does not imply the $\epsilon$ connectivity of $g^*$
(or $\delta$ connectivity for any $\delta$).  Therefore we do not have any estimate on the connectivity
of $g_\%\psi$ in (3) above.  It should be checked by an {\it ad hoc} method in each case.
For example, see \fullref{Exc}. 
\end{remarkun}

If the desuspension $\Omega C$ of a positive complex $C$ on $p_X$
is $\epsilon$ chain equivalent to a positive complex,
then $C$ is $\epsilon/4$ connected.  On the other hand, we have:

\begin{proposition}
\label{boundary}
\showlabel{boundary}
Let $n\ge 1$.

{\rm (1)}\qua Suppose an $n$--dimensional $\epsilon$
quadratic complex $(C,\psi)$ on $p_X$ is $\epsilon$ connected. Then $\d C$ is
$12\epsilon$ chain equivalent to an $(n{-}1)$--dimensional (resp. a 1--dimensional)
$4\epsilon$ chain complex $\hat\d C$ if $n > 1$ (resp. if\/ $n=1$).

{\rm (2)}\qua Suppose an $(n{+}1)$--dimensional $\epsilon$ quadratic pair
$(f\co C \to D, (\delta\psi, \psi))$ is $\epsilon$ connected.  Then
$\d D$ is $24\epsilon$ chain equivalent to an $n$--dimensional $5\epsilon$
chain complex $\hat\d D$.

{\rm (3)}\qua When $n=1$, the free 1--dimensional chain complex $(\hat\d C,1)$ given
in (1) and (2), viewed as a projective chain complex, is
$32\epsilon$ chain equivalent to a $0$--dimensional $32\epsilon$ projective
chain complex $(\tilde\d C,p)$ and there is a $32\epsilon$ isomorphism
\[
(\hat\d C_1,1) \oplus (\tilde\d C_0,p) \lga (\hat\d C_0,1),
\]
and hence the controlled reduced projective class $\bigl[\wtilde\d C,\!p\bigr]$ vanishes
in $\wwtilde K_0^{0,32\epsilon}\!(X;\!p_X,R)$.
\end{proposition}

\begin{proof} (1)\qua There exists a $4\epsilon$ morphism $h\co \d C_{-1} \to
\d C_0$ such that $dh \sim_{8\epsilon} 1$.
Define a $4\epsilon$ morphism $h'\co \d C_{n-1} \to \d C_n$ by
the composite:
\[
h'\co \d C_{n-1}=C_n\oplus C^1 \xrightarrow{\left(\begin{smallmatrix} 0 & (-)^n \\
(-)^n & 0\end{smallmatrix}\right)} C^1\oplus C_n \xrightarrow{~~h^*~~}  C^0=\d C_n,
\]
then $h'd \sim_{8\epsilon} 1$.
Now one can use the folding argument from the bottom (see Yamasaki
\cite{MY87}) using $h$ and, if $n>1$,
from the top (see Ranicki--Yamasaki \cite{RY-K}) using $h'$ to construct a desired chain equivalence.

(2)\qua There exists a $4\epsilon$ morphism $h\co \d D_{-1} \to \d D_0$
such that $dh \sim_{8\epsilon} 1$.
Define a $5\epsilon$ morphism $h'\co \d D_{n} \to \d D_{n+1}$ by
the composite of
\begin{equation*}
{\left(\begin{smallmatrix}0 & (-)^{n+1}&0 \\ (-)^{n+1} & 0 & 0\\ 0&0&(1+T)\psi_0
\end{smallmatrix}\right)}
\co \d D_{n}=D_{n+1}\oplus D^1\oplus C^0  \to
\d D^0=D^1\oplus D_{n+1}\oplus C_n
\end{equation*}
and $h^*\co \d D^0\to \d D^{-1}=\d D_{n+1}$,
then $h'd \sim_{8\epsilon} 1$.
Use the folding argument again.

(3)\qua The boundary map $\hat\d C_1=\d C_1\oplus \d C_{-1} \to \hat\d C_0 =\d C_0$
is given by the matrix $(d_{\d C} \quad h)$.  Therefore
\[
s={\begin{pmatrix} h' -h'hd_{\d C} \\ d_{\d C}\end{pmatrix}}
 \co \d C_0 \lga \d C_1 \oplus \d C_{-1}
\]
defines a $12\epsilon$ morphism $s\co \hat\d C_0 \to \hat\d C_1$ such that
$sd_{\hat\d C} \sim_{16\epsilon} 1$.  Define $\tilde\d C_0$ to be $\hat\d C_0$
and define a $16\epsilon$ morphism $p\co \tilde\d C_0 \to \tilde\d C_0$
by $1-d_{\hat\d C}$, then $p^2 \sim_{32\epsilon} p$ and
$p\co (\tilde\d C_0,1) \to (\tilde\d C_0,p)$ defines the desired $32\epsilon$
chain equivalence.
The isomorphism can be obtained by combining the following isomorphisms:
\begin{gather*}
\xymatrix{
(\hat\d C_1,1) \ar@<1ex>[r]^-{d} & (\hat\d C_0, 1-p) \ar@<1ex>[l]^-{s} }\\
\xymatrix{
(\hat\d C_0,1-p)\oplus (\hat\d C_0,p) \ar@<1ex>[r]^-{(q~~p)} & (\hat\d C_0, 1)
\ar@<1ex>[l]^-{{}^t(q~~p)}}
\end{gather*}
This completes the proof.
\end{proof}

Controlled connectivity is preserved under
union operation in the following manner.

\begin{proposition}
\label{union2}
\showlabel{union2}
If adjoining $\epsilon$ cobordisms $c$ and $c'$
are $\epsilon$ connected, then $c\cup c'$ is $100\epsilon$ connected.
\end{proposition}

\begin{proof} Similar to \fullref{union}.
\end{proof}

Now we define the epsilon-controlled $L$--groups.  Let $Y$ be a subset of $X$.

\begin{definition} For an integer $n\ge0$, pair of non-negative numbers $\delta\ge\epsilon\ge0$,
and a ring with involution $R$,
$L_n^{\delta,\epsilon}(X,Y;p_X,R)$ is defined to be the equivalence classes of
finitely generated $n$--dimensional $\epsilon$ connected $\epsilon$ quadratic complexes on $p_X$
that are $\epsilon$ Poincar\'e over $X-Y$.
The equivalence relation is generated by finitely generated $\delta$ connected $\delta$
cobordisms that are $\delta$ Poincar\'e over $X-Y$.
\end{definition}

\begin{remarkun} We use the following abbreviations for simplicity:
\begin{itemize}
\item $L_n^{\delta,\epsilon}(X;p_X,R) = L_n^{\delta,\epsilon}(X,\emptyset;p_X,R)$
\item $L_n^{\epsilon}(X,Y;p_X,R) = L_n^{\epsilon,\epsilon}(X,Y;p_X,R)$
\item $L_n^{\epsilon}(X;p_X,R) = L_n^{\epsilon,\epsilon}(X;p_X,R)$
\end{itemize}
\end{remarkun}

\begin{proposition}
\label{g-str}
\showlabel{g-str}
Direct sum $(C,\psi)\oplus (C',\psi')=(C\oplus C',\psi\oplus\psi')$
induces an abelian group structure on $L_n^{\delta,\epsilon}(X,Y;p_X,R)$.
Furthermore if $[C,\psi] = [C',\psi'] \in L_n^{\delta,\epsilon}(X,Y;p_X,R)$, then
there is a finitely generated $100\delta$ connected $2\delta$ cobordism between $(C,\psi)$
and $(C',\psi')$ that is $100\delta$ Poincar\'e over $X-Y^{100\delta}$.
\end{proposition}

\begin{proof} The inverse of an element $[C,\psi]$ is given
by $[C,-\psi]$.
In fact, as in \fullref{i-str} and \fullref{conn} (with $g=1$),
\[
((1 ~~1) \co C\oplus C \sa C, (0,\psi\oplus -\psi))
\]
gives an $\epsilon$ connected $\epsilon$ null-cobordism
of $(C,\psi)\oplus (C,-\psi)$
that is $\epsilon$ Poincar\'e over $X-Y$.
The second claim follows from \fullref{union} and \fullref{union2}, because we can glue a sequence
of cobordisms at once.
\end{proof}

If $\delta'\ge\delta$ and $\epsilon'\ge\epsilon$, then
there is a homomorphism
\[
L_n^{\delta,\epsilon}(X,Y;p_X,R) \lga L_n^{\delta',\epsilon'}(X,Y;p_X,R)
\]
which sends $[C,\psi]$ to $[C,\psi]$.
This is called the \emph{relax-control map}.

In the study of controlled $L$--groups, we need an analogue of \fullref{i-str} for pairs:

\begin{proposition}
\label{i-str2}
\showlabel{i-str2}
Suppose there is a triad of $\epsilon$ chain
complexes on $p_X$
\[
\xymatrix@R=5pt{
C \ar[dd]_{g} \ar[r]^{f} \ar @{~>} [ddr]^{g} &  D \ar[dd]^{h} \\
& & dk + kd\sim f'g - hf\\
C' \ar[r]_{f'} & D'
}
\]
where $f$, $f'$, $g$, $h$ are $\epsilon$ chain maps and $k$ is an $\epsilon$
chain homotopy, and suppose $(\delta\psi,\psi)$ is an $(n{+}1)$--dimensional
$\epsilon$ quadratic structure on $f$.

{\rm(1)}\qua There is induced a $4\epsilon$ quadratic structure on $f'$:
\begin{gather*}
(g,h;k)_\%(\delta\psi,\psi)=(h\delta\psi_sh^*+(-)^{n+1}k\psi_sf^*h^* +
    (-)^{n-r+1}f'g\psi_sk^*\\
+(-)^{r+1}kT\psi_{s+1}k^*\co (D'^{n+1-r-s},q'{}^*)\to (D'_r,q'_r),g\psi_sg^*)
    _{s\ge 0}.\\
\end{gather*}
{\rm(2)}\qua Suppose $g$ and $h$ are $\epsilon$ chain equivalences.
\begin{enumerate}
\item[\rm(a)] If $(\delta\psi,\psi)$ is $\epsilon$ Poincar\'e over $X-Y$,
then $(g,h;k)_\%(\delta\psi,\psi)$ is $30\epsilon$ Poincar\'e
over $X-Y^{81\epsilon}$.
\item[\rm(b)] If $(f, (\delta\psi,\psi))$ is $\epsilon$ connected,
then $(f',(g,h;k)_\%(\delta\psi,\psi))$ is $30\epsilon$ connected.
\end{enumerate}
\end{proposition}

\begin{proof} (1) is easy to check.  (2) can be checked by showing
that
\begin{align*}
((-)^{n+1-r}k\psi_0k^*\quad k(1+T)\psi_0g^*)\co
&\C(f')^{n+1-r}=D'{}^{n+1-r} \oplus C'{}^{n-r}\\
&\qquad \lga D'_{r+1}
\end{align*}
is a $3\epsilon$ chain homotopy between the duality map for
$(g,h;k)_\%(\delta\psi,\psi)$ and the following chain map:
\[
h((1+T)\delta\psi_0\quad f(1+T)\psi_0)
\begin{pmatrix}h^* &0\\ (-)^{n+1-r}k^*&g^*\end{pmatrix}
\co \C(f')^{n+1-r}\lga (D'_r,q'_r),
\]
which is a weak $27\epsilon$ chain equivalence over $X-Y^{18\epsilon}$
in case (2a), and $16\epsilon$ connected in case (2b).
\end{proof}

\begin{corollary}
\label{i-str3}
\showlabel{i-str3}
Suppose $f\co C\to D$ is an $\epsilon$ chain map,
$(\delta\psi,\psi)$ is an $(n{+}1)$--di\-men\-sion\-al $\epsilon$ quadratic structure
on $f$,  $g\co C\to C'$ is a $\gamma$ chain equivalence,
and $h\co D\to D'$ is a $\delta$ chain equivalence.
Let $\epsilon'=\gamma+\delta+\epsilon$ and $g^{-1}$ be a $\gamma$ chain
homotopy inverse of $g$.

{\rm(1)}\qua There is an $(n{+}1)$--dimensional $4\epsilon'$ quadratic structure
$(\delta\psi',\psi'=g_\%\psi)$ on the $\epsilon'$ chain map
$f'=hfg^{-1}\co (C',p')\to(D',q')$.

{\rm(2)}\qua If $(\delta\psi,\psi)$ is $\epsilon$ Poincar\'e over $X-Y$, then
$(\delta\psi', \psi')$ is $30\epsilon'$ Poincar\'e over $X-Y^{81\epsilon'}$.

(3)\qua If $(\delta\psi,\psi)$ is $\epsilon$ connected, then
$(\delta\psi', \psi')$ is $30\epsilon'$ connected.
\end{corollary}

\begin{proof}
Let $\Gamma\co g^{-1}g\simeq 1$ be a $\gamma$ chain homotopy.
Define an $\epsilon'$ chain homotopy $k\co hf\simeq f'g$ by $k=-hf\Gamma$,
and apply \fullref{i-str2}
\end{proof}

The last topic of this section is the functoriality.
A \emph{map between control maps} $p_X\co M\to X$ and $p_Y\co N\to Y$ means
a pair of continuous maps $(f\co M\to N, \wbar f\co X\to Y)$ which makes
the following diagram commute:
\[
\xymatrix{ M \ar[r]^f \ar[d]_{p_X}& N \ar[d]^{p_Y}\\
           X \ar[r]_{\wbar f}      & Y.}
\]
For example, given a control map $p_Y\co N\to Y$ and a subset $X\subset Y$,
let us denote the control map $p_Y|p_Y^{-1}(X) \co  p_Y^{-1}(X)\to X$
by $p_X\co M\to X$.  Then the inclusion maps $j\co M\to N$, $\wbar
\jmath\co X \to Y$ form a
map form $p_X$ to $p_Y$.

Epsilon controlled $L$--groups are functorial with respect to
maps and relaxation of control in the following sense.

\begin{proposition}
\label{func}
\showlabel{func}
Let $F=(f,\wbar f)$ be a map from
$p_X\co M\to X$ to $p_Y\co N\to Y$, and suppose that $\wbar f$ is Lipschitz
continuous with Lipschitz constant $\lambda$, ie
there exists a constant $\lambda>0$ such that
\[
d(\wbar f(x_1), \wbar f(x_2)) \le \lambda d(x_1,x_2)\qquad (x_1, x_2\in X).
\]
Then $F$ induces a homomorphism
\[
F_*\co L_n^{\delta,\epsilon}(X,X';p_X,R) \lga L_n^{\delta',\epsilon'}(Y,Y';p_Y,R)
\]
if $\delta' \ge \lambda\delta$, $\epsilon'\ge\lambda\epsilon$ and $\wbar f(X')\subset Y'$.
If two maps $F=(f,\wbar f)$ and $G=(g,\wbar g)$ are homotopic
through maps $H_t=(h_t,\wbar h_t)$ such that each $\wbar h_t$
is Lipschitz continuous with Lipschitz constant $\lambda$,
$\delta' > \lambda\delta$, $\epsilon' \ge \lambda\epsilon$,
and $\wbar h_t(X')\subset Y'$, then
$F$ and $G$ induce the same homomorphism:
\[
F_*=G_* \co L_n^{\delta,\epsilon}(X,X';p_X,R) \lga L_n^{\delta',\epsilon'}(Y,Y';p_Y,R).
\]

\end{proposition}

\begin{proof}
The direct image construction for geometric modules
and morphisms (see Ranicki--Yamasaki \cite[page 7]{RY-K}) can be used to define the direct images
$f_\#(C,\psi)$ of quadratic complexes and the direct images of cobordisms.
And this induces the desired $F_*$.

For the second part,
split the homotopy into thin layers to construct small cobordisms.
The size of the cobordism may be slightly bigger than the size of the
object itself.
\end{proof}

\begin{remarkun}
The above is stated for Lipschitz continuous maps to simplify the statement.
For specific $\delta\ge\epsilon$ and  $\delta'\ge\epsilon'$,
the following condition, instead of the Lipschitz condition above,
is sufficient for the existence of $F_*$:
\begin{align*}
&d(\wbar f(x_1),\wbar f(x_2))\le k\epsilon'\quad\hbox{~whenever~~~~}
d(x_1,x_2) \le k\epsilon, \hbox{ and}\\
&d(\wbar f(x_1),\wbar f(x_2))\le k\delta'\quad\hbox{~whenever~~~~}
d(x_1,x_2) \le k\delta,
\end{align*}
for $k = 1$, 3, 4, 8.
The second part of the proposition also holds under this condition.
When $X$ is compact and $\delta'\ge\epsilon'$ are given,
the uniform continuity of $\wbar f$ implies that this condition is satisfied
for sufficiently small pairs $\delta\ge\epsilon$.
\end{remarkun}

\section{Epsilon-controlled projective $L$--groups}\label{proj-L}

\showlabel{proj-L}
Fix a subset $Y$ of $X$, and let $\F$ be a family of subsets of $X$ such that
$Z\supset Y$ for each $Z\in\F$.
In this section we introduce intermediate epsilon-controlled $L$--groups
$L_n^{\F,\delta,\epsilon}(Y;p_X,R)$, which will appear in the stable-exact
sequence of a pair (\fullref{relative-L}) and
also in the Mayer--Vietoris sequence (\fullref{MV}).
Roughly speaking, these are defined using ``controlled projective quadratic
chain complexes'' $((C,p),\psi)$ with vanishing
$\epsilon$--controlled reduced projective class
$[C,p]=0\in \wtilde K_0^{n,\epsilon}(Z;p_Z,R)$ for each $Z\in\F$.

$\wtilde K_0^{n,\epsilon}(Z;p_Z,R)$ is an abelian group defined
as the set of equivalence classes  $[C,p]$ of finitely generated
$\epsilon$ projective chain complexes on $p_Z$. See Ranicki--Yamasaki \cite{RY-K}
for the details. The following is known \cite[3.1 and 3.5]{RY-K}:

\begin{proposition}
\label{zero-in-K}
\showlabel{zero-in-K}
If $[C,p]=0\in\wtilde K_0^{n,\epsilon}(Z;p_Z,R)$, then there is an
$n$--dimensional free $\epsilon$ chain complex $(E,1)$ such that
$(C,p)\oplus (E,1)$ is $3\epsilon$ chain equivalent to an
$n$--dimensional free $\epsilon$ chain complex on $p_Z$.
If we further assume that $n\ge 1$, then $(C,p)$ itself is $60\epsilon$ chain
equivalent to an $n$--dimensional free $30\epsilon$ chain complex on $p_Z$.
\end{proposition}

All the materials in the previous two sections (except for \fullref{boundary}(3)) have obvious
analogues in the category of projective chain complexes with the identity maps
in the formulae replaced by appropriate projections.  So we shall only describe
the basic definitions and omit stating the obvious analogues of
Propositions~\ref{comp}--\ref{cont-triad}, \fullref{P-P},
Propositions~\ref{union} and~\ref{i-str}, Lemmas~\ref{l-cont}
and~\ref{conn}, Propositions~\ref{union2} and~\ref{i-str2}, and
\fullref{i-str3},
and we refer them by \fullref{comp}${}'$, \fullref{cone}${}'$, \dots.
An analogue of \fullref{g-str} will be explicitly stated in
\fullref{g-str2} below.

For a projective module $(A,p)$ on $p_X$, its dual $(A,p)^*$ is the
projective module $(A^*,p^*)$ on $p_X$.
If $f\co (A,p)\to (B,q)$ is an $\epsilon$ morphism \cite{RY-K},
then $f^*\co (B,q)^*\to(A,p)^*$ is also an $\epsilon$ morphism.
For an $\epsilon$ projective chain complex on $p_X$
\[
(C,p)\co  \dots \sa(C_r,p_r)\xrightarrow{\ d_r\ } (C_{r-1},p_{r-1})\xrightarrow{d_{r-1}}\dots
\]
in the sense of \cite{RY-K}, $(C,p)^{n-*}$ will denote the $\epsilon$ projective
chain complex on $p_X$ defined by:
\[
\dots\sa(C^{n-r},p_{n-r}^*)\xrightarrow{(-)^rd_r^*}(C^{n-r+1},p_{n-r+1}^*)\sa\dots.
\]

Before we go on to define $\epsilon$ projective quadratic complexes, we need
to define basic notions for projective chain complexes.  For $Y=X$ or
for free chain complexes, these are already defined in \cite{RY-K}.

Suppose $f\co (A,p)\to (B,q)$ is a morphism between projective
modules on $p_X$, and let $Y$ be a subset of $X$.  The {\sl
restriction} $f|Y$ of $f$ to $Y$ will mean the restriction of $f$
in the sense of \cite[page 21]{RY-K} with $f$ viewed as a geometric
morphism from $A$ to $B$; that is, $f|Y$ is the sum of the
paths (with coefficients) that start from points in $p_X^{-1}(Y)$.
$f|Y$ can be viewed as a geometric morphism from $A$ to $B$ and
also as a geometric morphism from $A(Y)$ to $B(Y)$, where $A(Y)$
denotes the restriction of $A$ to $Y$ in the sense of \cite{RY-K},
ie the geometric submodule of $A$ generated by the basis
elements of $A$ that are in $p_X^{-1}(Y)$. But, in general, it
does not give a morphism from $(A,p)$ to $(B,q)$. Also note that
there is no obvious way to ``restrict'' a projection $p\co A\to A$ to
a projection on $A(Y)$.

The following four paragraphs are almost verbatim copies of the
definitions for free chain complexes \cite[page 22]{RY-K}.

Let $f$, $g\co (A,p)\to (B,q)$ be morphisms; $f$ is said to be
\emph{equal to $g$ over $Y$} ($f=g$ \emph{over} $Y$) if $f|Y= g|Y$, and
$f$ is said to be \emph{$\epsilon$ homotopic to $g$ over $Y$}
($f\sim_{\epsilon}g$ \emph{over} $Y$)
if $f|Y\sim_\epsilon g|Y$.

Let $f$, $g\co (C,p)\to (D,q)$ be chain maps between projective chain complexes on
$p_X$.  A collection $\{h_r\co (C_r,p_r)\to (D_{r+1},q_{r+1})\}$ of $\epsilon$
morphisms is said to be an \emph{$\epsilon$ chain homotopy over $Y$ between
$f$ and $g$} if $dh+hd\sim_{2\epsilon} g-f$ over $Y$.

An $\epsilon$ chain map $f\co (C,p)\to (D,q)$ is said to be an \emph{$\epsilon$
chain equivalence over $Y$} if there exist an $\epsilon$ chain map $g\co (D,q)
\to (C,p)$ and $\epsilon$ chain homotopies over $Y$ between $gf$ and $p$ and
between $fg$ and $q$.

A chain complex $(C,p)$ is said to be \emph{$\epsilon$ contractible over $Y$}
if there is an $\epsilon$ chain homotopy over $Y$ between $0\co (C,p)\to (C,p)$
and $p\co (C,p)\to (C,p)$; such a chain homotopy over $Y$ is called an
\emph{$\epsilon$ chain contraction of $(C,p)$ over $Y$}.

The \fullref{weak-eq} of \emph{weak $\epsilon$ chain equivalences over $Y$}
(for chain maps between free chain complexes) can be
rewritten for maps between projective chain complexes in the obvious manner.

The following is the most important technical proposition
in the theory of controlled projective chain complexes.

\begin{proposition}[Ranicki--Yamasaki {\cite[5.1 and 5.2]{RY-K}}]
\label{5.1-5.2}
\showlabel{5.1-5.2}
If an $n$--dimensional free $\epsilon$ chain complex
$C$ on $p_X$ is $\epsilon$ contractible over $X-Y$, then $(C,1)$ is
$(6n+15)\epsilon$ chain equivalent to an $n$--dimensional $(3n+12)\epsilon$
projective chain complex on $p_{Y^{(4n+14)\epsilon}}$.
Conversely, if an $n$--dimensional free chain complex $(C,1)$ on $p_X$ is
$\epsilon$ chain equivalent to a projective chain complex $(D,r)$ on
$p_Y$, then $C$ is $\epsilon$ contractible over $X-Y^\epsilon$.
\end{proposition}

Now we introduce quadratic structures on projective
chain complexes and pairs. An \emph{$n$--dimensional $\epsilon$
quadratic structure} on a projective chain complex $(C,p)$ on
$p_X$ is an $n$--dimensional $\epsilon$ quadratic structure $\psi$
on $C$ (in the sense of \fullref{q-str}) such that
$\psi_s\co (C^{n-r-s},p^*)\to (C_r,p)$ is an $\epsilon$ morphism for
every $s\ge0$ and $r\in\ZZ$. Similarly, an {\sl
$(n{+}1)$--dimensional $\epsilon$ quadratic structure} on a chain map
$f\co (C,p)\to (D,q)$ is an $(n{+}1)$--dimensional $\epsilon$ quadratic
structure $(\delta\psi,\psi)$ on $f\co C\to D$ such that
$\delta\psi_s\co (D^{n+1-r-s},q^*) \to (D_r,q)$ and
$\psi_s\co (C^{n-r-s},p^*)\to (C_r,p)$ are $\epsilon$ morphisms for
every $s\ge0$ and $r\in\ZZ$. An $n$--dimensional $\epsilon$
projective chain complex $(C,p)$ on $p_X$ equipped with an
$n$--dimensional $\epsilon$ quadratic structure is called an {\sl
$n$--dimensional $\epsilon$ projective quadratic complex on} $p_X$,
and an $\epsilon$ chain map $f\co (C,p)\to (D,q)$ between an
$n$--dimensional $\epsilon$ projective chain complex $(C,p)$ on
$p_X$ and an $(n{+}1)$--dimensional $\epsilon$ projective chain
complex $(D,q)$ on $p_X$ equipped with an $(n{+}1)$--dimensional
$\epsilon$ quadratic structure is called an {\sl
$(n{+}1)$--dimensional $\epsilon$ projective quadratic pair on}
$p_X$.

An \emph{$\epsilon$ cobordism} of $n$--dimensional $\epsilon$ projective
quadratic complexes $((C,p),\psi)$, $((C',p'),\allowbreak\psi')$ on $p_X$ is
an  $(n{+}1)$--dimensional $\epsilon$ projective quadratic pair on $p_X$
\[
((f ~~f')\co (C,p)\oplus(C',p')\sa (D,q),(\delta\psi,\psi\oplus-\psi'))
\]
with boundary $((C,p)\oplus(C',p'),\psi\oplus-\psi')$.

Boundary constructions, algebraic Poincar\'e thickenings, algebraic Thom
complexes, $\epsilon$ connectedness are defined as in the previous section.

An $n$--dimensional $\epsilon$ quadratic structure $\psi$ on $(C,p)$ is
$\epsilon$ \emph{Poincar\'e (over $Y$)} if
\[
\d(C,p)=\Omega\C((1+T)\psi_0\co (C^{n-*},p^*)\sa (C,p))
\]
is $4\epsilon$ contractible (over $Y$).
$((C,p),\psi)$ is $\epsilon$ \emph{Poincar\'e (over $Y$)}
if $\psi$ is $\epsilon$ Poincar\'e (over $Y$).
Similarly, an $(n{+}1)$--dimensional $\epsilon$ quadratic structure $(\delta\psi,
\psi)$ on $f\co (C,p)\to (D,q)$ is $\epsilon$ \emph{Poincar\'e (over $Y$)}
if $\d(C,p)$ and
\[
\d(D,q)=\Omega\C(((1+T)\delta\psi_0\quad f(1+T)\psi_0)\co \C(f)^{n+1-*}\sa (D,q))
\]
are both $4\epsilon$ contractible (over $Y$).
A pair $(f,(\delta\psi,\psi))$ is $\epsilon$ \emph{Poincar\'e (over $Y$)}
if $(\delta\psi,\psi)$ is $\epsilon$ Poincar\'e (over $Y$).

Let $Y$ and be a subset of $X$ and $\F$ be a family of subsets of $X$ such that
$Z\supset Y$ for every $Z\in\F$.

\begin{definition}
\label{Lp}
\showlabel{Lp}
 Let $n\ge 0$ and $\delta\ge\epsilon\ge0$.
$L_n^{\F,\delta,\epsilon}(Y;p_X,R)$ is the equivalence
classes of finitely generated $n$--dimensional $\epsilon$ Poincar\'e
$\epsilon$ projective quadratic complexes $((C,p),\psi)$ on $p_Y$
such that $[C,p]= 0$ in $\wtilde K_0^{n,\epsilon}(Z;p_Z,R)$
for each $Z\in\F$.
The equivalence relation is generated by finitely generated
$\delta$ Poincar\'e $\delta$ cobordisms $((f~~f')\co (C,p)\oplus (C',p') \to (D,q),
(\delta\psi,\psi\oplus -\psi'))$ on $p_Y$
such that $[D,q]=0$ in $\wtilde K_0^{n+1,\delta}(Z;p_Z,R)$ for each $Z\in\F$.
\end{definition}

\begin{remarkun} We use the following abbreviation:
$L_n^{\F, \epsilon}(Y;p_X,R) = L_n^{\F,\epsilon,\epsilon}(Y;p_X,R)$.
\end{remarkun}

\begin{proposition}
\label{g-str2}
\showlabel{g-str2}
Direct sum induces an abelian group structure
on $L_n^{\F,\delta,\epsilon}(Y;p_X,R)$.
Furthermore if $[(C,p),\psi] = [(C',p'),\psi'] \in L_n^{\F,\delta,\epsilon}(Y;p_X,R)$,
then there is a finitely generated $100\delta$ Poincar\'e $2\delta$ cobordism on $p_Y$
\[
((f~~f')\co (C,p)\oplus (C',p') \to (D,q), (\delta\psi,\psi\oplus -\psi'))
\]
such that $[D,q]=0$ in $\wtilde K_0^{n+1,9\delta}(Z;p_Z,R)$ for each $Z\in\F$.
\end{proposition}

\begin{proof}
The first part is similar to the proof of \fullref{g-str}. Observe that
$[D,q]=0$ in $\wtilde K_0^{n+1,9\delta}(Z;p_Z,R)$, because
$$[\C(g\co (E,r)\to(F,s))] = [F,s]-[E,r]\in \wtilde
K_0^{n+1,9\delta}(Z;p_Z,R)$$
for any $\delta$ chain map $g$ between $\delta$ projective chain complexes
$(E,r)$ (of dimension $n$) and $(F,s)$ (of dimension $n+1$) on $p_Z$.
See Ranicki--Yamasaki \cite[page 18]{RY-K}.
\end{proof}

A functoriality with respect to maps and relaxation of control similar to
\fullref{func} holds for epsilon-controlled projective $L$--groups.

\begin{proposition}
\label{func2}
\showlabel{func2}
Let $F=(f,\wbar f)$ be a map from
$p_X\co M\to X$ to $p_Y\co N\to Y$, and suppose that $\wbar f$ is Lipschitz
continuous with Lipschitz constant $\lambda$, ie
there exists a constant $\lambda>0$ such that
\[
d(\wbar f(x_1), \wbar f(x_2)) \le \lambda d(x_1,x_2)\qquad (x_1, x_2\in X).
\]
If $\delta' \ge \lambda\delta$, $\epsilon' \ge \lambda\epsilon$, $\wbar f(A)\subset B$,
and there exists a $Z\in\F$ satisfying $\wbar f(Z)\subset Z'$
for each $Z'\in \F'$,  then $F$ induces a homomorphism
\[
F_*\co L_n^{\F,\delta,\epsilon}(A;p_X,R) \lga L_n^{\F',\delta',\epsilon'}(B;p_Y,R).
\]
It is $\lambda$--Lipschitz-homotopy invariant if $\delta'>\lambda\delta$ in addition.
\end{proposition}

\begin{remarkun}
As in the remark to \fullref{func}, for a specific $\delta$ and  $\epsilon$,
we do not need the full Lipschitz condition to guarantee the existence of $F_*$.
\end{remarkun}

There is an obvious homomorphism
\[
\iota\co L_n^{\delta,\epsilon}(Y;p_Y,R) \lga L_n^{\F,\delta,\epsilon}(Y;p_X,R);
    \qquad [C,\psi] \mapsto [(C,1), \psi]
\]
from free $L$--groups to projective $L$--groups.
On the other hand,
the controlled $K$--theoretic condition posed in the definition
can be used to construct homomorphisms from  projective $L$--groups
to free $L$--groups:

\begin{proposition}
\label{P2F}
\showlabel{P2F}
There exist a constant $\alpha > 1$
such that the following holds true:
for any control map $p_X\co M\to X$, any subset $Y\subset X$, any family of
subsets $\F$ of $X$ containing $Y$ , any element $Z$ of $\F$,
any number $n\ge 0$, and any pair of positive numbers $\delta\ge\epsilon$
and $\delta\ge\epsilon$ with $\delta'\ge\alpha\delta$, $\epsilon'\ge\alpha\epsilon$,
there is a well-defined homomorphism
\[
\iota_Z\co L_n^{\F,\delta,\epsilon}(Y;p_X,R) \lga
  L_n^{\delta',\epsilon'}(Z;p_Z,R),
\]
functorial with respect to relaxation of control,
such that the composite maps
\begin{align*}
&L_n^{\F,\delta,\epsilon}(Y;p_X,R)\xrightarrow{\iota_Z}
L_n^{\delta',\epsilon'}(Z;p_Z,R)\xrightarrow{\iota}
L_n^{p,\delta',\epsilon'}(Z;p_Z,R)\\
&L_n^{\delta,\epsilon}(Y;p_Y,R)\xrightarrow{\iota}
L_n^{\F,\delta,\epsilon}(Y;p_X,R)\xrightarrow{\iota_Z}
L_n^{\delta',\epsilon'}(Z;p_Z,R)
\end{align*}
are equal to the ones induced from inclusion maps.
\end{proposition}

\begin{remarkun}
Actually $\alpha=20000$ works.
In the rest of the paper, we always assume that $\alpha=20000$.
\end{remarkun}

\begin{proof}
Let $[(C,p),\psi]$ be an element of $L_n^{\F,\delta,\epsilon}(Y;p_X,R)$,
and fix $Z\in\F$.
Recall that $[C,p]=0 \in \wtilde K_0^{n,\epsilon}(Z;p_Z,R)$.
By  \fullref{zero-in-K}, there exists an $n$--dimensional free
$\epsilon$ chain complex $(E,1)$ on $p_Z$ such that
$(C,p)\oplus (E,1)$ is $3\epsilon$ chain equivalent to some $n$--dimensional
free $\epsilon$ chain complex $(\wbar F,1)$ on $p_Z$.
Add $1\co (E^{n-*},1) \to (E^{n-*},1)$ to this chain equivalence to get a
$3\epsilon$ chain equivalence
\[
g\co (C,p)\oplus (\d\Sigma E,1) \lga (\wbar F,1)\oplus (E^{n-*},1)=(F,1)
\]
of projective chain complexes on $p_Z$,
where $\Sigma E$ is defined using the trivial $(n{+}1)$--dimensional quadratic structure
$\theta=0$ on $\Sigma E$. See \fullref{b-ex}.
We set
\[
\iota_Z [(C,p),\psi]=[F,g_\%(\psi\oplus\d\theta)].
\]
Let us show that this defines a well-defined map.
Suppose $[(C,p),\psi]=[(C',p'),\psi']$ in $L_n^{\F,\delta,\epsilon}(Y;p_X,R)$,
and let $E$ and $E'$ be $n$--dimensional free $\epsilon$ chain complexes on $p_Z$
together with $3\epsilon$ chain equivalences
\begin{gather*}
g\co (\wbar C,\wbar p)=(C,p)\oplus(\d\Sigma E,1) \to (F,1)\\
g'\co (\wbar C',\wbar p')=(C',p')\oplus(\d\Sigma E',1) \to (F',1)
\end{gather*}
to free $\epsilon$ chain complexes $F$ and $F'$ on $p_Z$.
By \fullref{g-str2} above and \fullref{zero-in-K},
there is a $100\delta$ Poincar\'e $2\delta$
null-cobordism
\[
(f\co (C,p)\oplus(C',p') \sa (D,q), (\delta\psi,\psi\oplus -\psi'))
\]
such that $(D,q)$ is $540\delta$ chain equivalent to an $(n{+}1)$--dimensional
free $270\delta$ chain complex $(G,1)$ (as a projective chain complex on $p_Z$).
Take the direct sum with the null-cobordisms
\begin{gather*}
(i_{\Sigma E}\co (\d\Sigma E,1) \sa (E^{n-*},1), (0,\d\theta)),\\
(i_{\Sigma E'}\co (\d\Sigma E',1) \sa (E'^{n-*},1), (0,-\d\theta')).
\end{gather*}
Now the claim follows from \fullref{i-str3}${}'$.
\[
\xymatrix@C+30pt{
(\wbar C,\wbar p)\oplus (\wbar C',\wbar p')\ar[r]^-{100\delta \text{ Poincar\'e}}
\ar[d]_-{g\oplus g'}^-{\simeq_{3\epsilon}} &
(D,q)\oplus (E^{n-*},1)\oplus(E'^{n-*},1)\ar[d]^-{\simeq_{540\delta}}\\
(F,1)\oplus (F',1) \ar@{.>}[r]_-{20000\delta \text{ Poincar\'e}}
&(G,1)\oplus (E^{n-*},1)\oplus(E'^{n-*},1)}
\]
This completes the proof.
\end{proof}

\section{Stably-exact sequence of a pair}\label{relative-L}

\showlabel{relative-L}
Let $Y$ be a subset of $X$.
We discuss relations between the various controlled $L$--groups of
$X$, $Y$, and $(X,Y)$ by fitting them into a stably-exact sequence.
Two of the three kinds of maps forming the sequence have already appeared.
The first is the map
\[
i_*=\iota_X\co L_n^{\{X\},\delta,\epsilon}(Y;p_X,R) \lga L_n^{\delta',\epsilon'}(X;p_X,R)
\]
defined when $\delta'\ge\alpha\delta$ and $\epsilon'\ge\alpha\epsilon$.
The second is the homomorphism induced by the inclusion map $j\co (X,\emptyset)\to(X,Y)$:
\[
j_*\co L_n^{\delta,\epsilon}(X;p_X,R) \to L_n^{\delta',\epsilon'}(X,Y;p_X,R).
\]
defined for positive numbers $\delta'\ge\delta$ and $\epsilon'\ge\epsilon$.
The third map $\d$ is described in the next proposition.

\begin{proposition}
\label{w-def-d}
\showlabel{w-def-d}
For $n\ge 1$, there exists a constant $\kappa_n>1$ such that the following holds true:
If $Y'\supset Y^{\kappa_n\delta}$,
$\delta'\ge \kappa_n\delta$, and $\epsilon'\ge \kappa_n\epsilon$,
$\d([C,\psi])=[(E,q),\beta_\%\d\psi]$ defines a well-defined homomorphism:
\[
\d\co L_n^{\delta,\epsilon}(X,Y;p_X,R) \to L_{n-1}^{\{X\},\delta',\epsilon'}(Y';p_X,R),
\]
where
\[
\beta \co (\d C,1) \lga (E,q)
\]
is any $(200n+300)\epsilon$ chain equivalence
from $(\d C,1)$ to some $(n{-}1)$--dimensional $(100n+300)\epsilon$ projective
chain complex on $p_{Y'}$.
\end{proposition}

\begin{remarkun}
Actually $\kappa_n=150000(n+2)$ works.
In the rest of the paper, we always assume that $\kappa_n\ge 150000(n+2)$.
\end{remarkun}

\begin{proof}
We first show the existence of such $\beta$.
Take $[C,\psi]\in L_n^{\delta,\epsilon}(X,Y;p_X,R)$.
Suppose $n>1$.  By \fullref{boundary}(1),
there is a $12\epsilon$
chain equivalence between $\d C$ and an $(n{-}1)$--dimensional $4\epsilon$
chain complex $\hat \d C$ on $p_X$.  Since $\d C$ is $4\epsilon$ contractible
over $X-Y$, $\hat\d C$ is $28\epsilon$ contractible over $X-Y^{12\epsilon}$
by \fullref{l-cont}.
Now by \fullref{5.1-5.2}, $(\hat\d C,1)$ is $(168n+252)\epsilon$
($=(6(n-1)+15)\times 28\epsilon$)
chain equivalent
to an $(n{-}1)$--dimensional $(84n+252)\epsilon$ projective chain complex on
$p_{Y^{(112n+292)\epsilon}}$.

Next suppose  $n=1$.  By \fullref{boundary}(1) and (3), there is a $44\epsilon$
chain equivalence between $(\d C,1)$ and a $0$--dimensional $32\epsilon$
chain complex $(\tilde \d C,p)$.  Since $\d C$ is $4\epsilon$ contractible
over $X-Y$, $(\tilde\d C,p)$ is $92\epsilon$ contractible over
$X-Y^{44\epsilon}$, ie $p\sim_{184\epsilon}0$ over $X-Y^{44\epsilon}$.
Let $E=\tilde\d C | Y^{76\epsilon}$ and $q=p|Y^{44\epsilon}$, then
$p-q=p|(X-Y^{44\epsilon}) \sim_{184\epsilon}0$.  Therefore
\[
q~\sim_{184\epsilon}~p~\sim_{32\epsilon}~p^2~\sim_{216\epsilon}~pq
~\sim_{216\epsilon}~q^2,
\]
and $(E,q)$ is a 0--dimensional $216\epsilon$ projective chain complex on
$p_{Y^{292\epsilon}}$.
The $32\epsilon$ morphism $q$ defines a $216\epsilon$ isomorphism between
$(\tilde\d C,p)$ and $(E,q)$ in each direction.  Therefore $(\d C,1)$
is $260\epsilon$ chain equivalent to $(E,q)$.  This completes the
proof of the existence of $\beta$.

Suppose $[C,\psi]=[C'\psi']\in L_n^{\delta,\epsilon}(X,Y;p_X,R)$
and suppose $\beta\co (\d C,1) \to (E,q)$ and $\beta'\co (\d C',1) \to (E',q')$ are
chain equivalences satisfying the condition, and suppose $Y'$, $\delta'$, and $\epsilon'$
satisfy the hypothesis.
We show that $((E,q),\beta_\%\d\psi)$ and $((E',q'),\beta'_\%\d\psi')$
represent the same element in $L^{\{X\},\delta',\epsilon'}_{n-1}(Y';p_X,R)$.
Without loss of generality, we may assume that there is an $\epsilon$ connected
$\epsilon$ cobordism
\[
((f\quad f')\co C\oplus C' \sa D,\ (\d\psi,\psi\oplus -\psi'))
\]
which is $\epsilon$ Poincar\'e over $X-Y$. Apply the boundary
construction (\fullref{L-groups})to this pair to get a $3\epsilon$
Poincar\'e $2\epsilon$ quadratic structure $(\d \delta\psi,
\d\psi\oplus -\d\psi')$ on the $2\epsilon$ chain map $(\d
C,1)\oplus (\d C',1) \sa (\d D,1)$ of $2\epsilon$ chain complexes.
By \fullref{l-cont}${}'$, \fullref{boundary} and \fullref{5.1-5.2}, $(\d D,1)$
is $(312n+904)\epsilon$ chain equivalent to an $n$--dimensional
$(156n+624)\epsilon$ projective chain complex $(F,r)$ on
$p_{Y^{(208n+752)\epsilon}}$. Now, by \fullref{i-str3}${}'$, we can
obtain a $(15360n+36210)\epsilon$ Poincar\'e cobordism
\[
(E,q)\oplus(E',q')~ \lga~(F,r),\qquad (\chi, \beta_\%\d\psi\oplus(-\beta'_\%\d\psi')).
\]
Since such a structure involves $8(15360n+36210)\epsilon$ homotopies,
this cobordism can be viewed to be on $p_{Y^{(123088n+290432)\epsilon}}$.
Also $[F,r]=[\hat\d D,1]=0$ in $\wtilde K_0^{n,\epsilon'}(X;p_X,R)$,
and similarly $[E,q]=[E',q']=0$ in $\wtilde K_0^{n-1,\epsilon'}(X;p_X,R)$.
Therefore $[(E,q),\beta_\%\d\psi]=[(E',q'),\beta'_\%\d\psi']$ in
$L^{\{X\},\delta',\epsilon'}_{n-1}(Y';p_X,R)$.

\end{proof}

\begin{theorem}
\label{ThESP}\showlabel{ThESP} For any integer $n\ge 0$, there exists a constant
$\lambda_n > 1$ which depends only on $n$ such that the following holds true
for any control map $p_X$ and a subset $Y$ of $X$:

{\rm (1)}\qua Suppose $\delta'\ge \alpha\delta$, $\epsilon'\ge \alpha\epsilon$,
$\delta''\ge\delta'$, and $\epsilon''\ge\epsilon'$ so that
the following two maps are defined:
\[
L_n^{\{X\},\delta,\epsilon}(Y;p_X,R)\xrightarrow{~i_*~} L_n^{\delta',\epsilon'}(X;p_X,R)
\xrightarrow{~j_*~} L_n^{\delta'',\epsilon''}(X,W;p_X,R).
\]
If $W\supset Y^{\alpha\epsilon}$, then $j_*i_*$ is zero.

{\rm (2)}\qua
Suppose $\delta''\ge\delta'$, $\epsilon''\ge\epsilon'$ so that
$j_*\co L_n^{\delta',\epsilon'}(X;p_X,R) \to L_n^{\delta'',\epsilon''}(X,W;p_X,R)$ is defined.
If $\delta\ge\lambda_n\delta''$ and $Y\supset W^{\lambda_n\delta''}$, then the relax-control
image of the kernel of $j_*$ in $L_n^{\alpha\delta}(X;p_X,R)$ is
contained in the image of $i_*$ below:
\[
\xymatrix{
 & L_n^{\delta',\epsilon'}(X;p_X,R) \ar[r]^{j_*~~} \ar[d]
 & L_n^{\delta'',\epsilon''}(X,W;p_X,R) \\
L_n^{\{X\},\delta}(Y;p_X,R) \ar[r]^{i_*}
 & L_n^{\alpha\delta}(X;p_X,R) & \\
}
\]
{\rm (3)}\qua
Suppose $n\ge 1$, $\delta'\ge\delta$, $\epsilon'\ge\epsilon$,
$W\supset Y^{\kappa_n\delta'}$, $\delta''\ge\kappa_n\delta'$,
and $\epsilon''\ge\kappa_n\epsilon'$
so that the following two maps are defined:
\[
L_n^{\delta,\epsilon}(X;p_X,R)\xrightarrow{j_*}{}
L_n^{\delta',\epsilon'}(X,Y;p_X,R) \xrightarrow{\ \d\ }{}
L_{n-1}^{\{X\},\delta'',\epsilon''}(W;p_X,R).
\]
Then $\d j_*$ is zero.

{\rm (4)}\qua
Suppose $n\ge 1$, $W\supset Y^{\kappa_n\delta'}$, $\delta''\ge \kappa_n\delta'$,
and $\epsilon''\ge\kappa_n\epsilon'$ so that the map
$$\d\co L_n^{\delta',\epsilon'}(X,Y;p_X,R) \to L_{n-1}^{\{X\},
\delta'',\epsilon''}(W;p_X,R)$$
is defined.
If $\delta\ge\lambda_n\delta''$ and $Y'\supset W^{\lambda_n\delta''}$, then the relax-control
image of the kernel of $\d$ in $L_n^{\delta}(X,Y';p_X,R)$ is
contained in the image of $j_*$ below:
\[
\xymatrix{
 & L_n^{\delta',\epsilon'}(X,Y;p_X,R) \ar[r]^-{\d} \ar[d]
 & L_{n-1}^{\{X\}, \delta'',\epsilon''}(W;p_X,R) \\
L_n^{\delta}(X;p_X,R) \ar[r]^{j_*}
 & L_n^{\delta}(X,Y';p_X,R) & \\
}
\]
{\rm (5)}\qua
Suppose $n\ge 1$, $Y'\supset Y^{\kappa_n\delta}$, $\delta'\ge\kappa_n\delta$,
$\epsilon'\ge\kappa_n\epsilon$,
$\delta''\ge\alpha\delta'$, and $\epsilon''\ge\alpha\epsilon'$
so that the following two maps are defined:
\[
L_n^{\delta,\epsilon}(X,Y;p_X,R) \xrightarrow{\ \d\ }{}
L_{n-1}^{\{X\},\delta',\epsilon'}(Y';p_X,R)\xrightarrow{i_*}{}
L_{n-1}^{\delta'',\epsilon''}(X;p_X,R).
\]
Then $i_*\d$ is zero.

{\rm (6)}\qua
Suppose $n\ge 1$, $\delta''\ge \alpha\delta'$, and $\epsilon''\ge\alpha\epsilon'$ so that
$i_*\co L_{n-1}^{\{X\},\delta',\epsilon'}(Y;p_X,R) \to
L_{n-1}^{\delta'',\epsilon''}(X;p_X,R)$ is defined.
If $\delta\ge\lambda_n\delta''$ and $W\supset Y^{\lambda_n\delta''}$, then the relax-control
image of the kernel of $i_*$ in $L_{n-1}^{\{X\},\kappa_n\delta}(W^{\kappa_n\delta};p_X,R)$ is
contained in the image of $\d$ below:
\[
\xymatrix{
 & L_{n-1}^{\{X\},\delta',\epsilon'}(Y;p_X,R) \ar[r]^{i_*} \ar[d]
 & L_{n-1}^{\delta'',\epsilon''}(X;p_X,R) \\
L_n^{\delta}(X,W;p_X,R) \ar[r]^-{\d}
 & L_{n-1}^{\{X\}, \kappa_n\delta}(W^{\kappa_n\delta};p_X,R) & \\
}
\]
\end{theorem}

\begin{proof}
(1)\qua Let $[(C,p),\psi]\in L_n^{\{X\},\delta,\epsilon}(Y;p_X,R)$.
There is a $3\epsilon$ chain equivalence $g\co (C,p)\oplus(\d\Sigma E)\to (F,1)$
for some $n$--dimensional free $\epsilon$ chain complexes $E$ and $F$ on $p_X$,
and $j_*i_*[(C,p),\psi]\in L_n^{\delta'',\epsilon''}(X,W;p_X,R)$ is represented by
$(F,g_\%(\psi\oplus\d\theta))$, where $\theta$ is the trivial quadratic
structure on $\Sigma E$.
Take the sum of
\[
(0\co (C,p)\to 0,(0,\psi)), \quad{\rm and}\quad
(i_{\Sigma E}\co (\d\Sigma E,1) \to (E^{n-*},1), (0,\d\theta)).
\]
$(0,\psi\oplus\d\theta)$ is a $2\epsilon$ connected $2\epsilon$ quadratic
structure, and it is $2\epsilon$ Poincar\'e over $X-Y$.
Use the chain equivalence $g$ and \fullref{i-str3}${}'$ to get a $180\epsilon$ connected
$24\epsilon$ null-cobordism
\[
(F\lga E^{n-*},~~(\chi, g_\%(\psi\oplus\d\theta)))
\]
that is $180\epsilon$ Poincar\'e over $X-Y^{486\epsilon}$.

(2)\qua Let $[C,\psi]\in L_n^{\delta',\epsilon'}(X;p_X,R)$ and assume $j_*[C,\psi]=0\in
L_n^{\delta'',\epsilon''}(X,W;p_X,R)$.
By \fullref{g-str}, there is a $100\delta''$ connected $2\delta''$ null-cobordism
\[
(f\co C\to D,(\delta\psi,\psi))
\]
that is $100\delta''$ Poincar\'e over $X-W^{100\delta''}$.
Apply the boundary construction to this null-cobordism to get a $4\delta''$
chain map $\d f$ of $4\delta''$ chain complexes and an $n$--dimensional
$6\delta''$ Poincar\'e $6\delta''$ quadratic structure on it:
\[
\d f\co \d C \to \d D, \qquad \Psi_3=(\d\delta\psi,\d\psi).
\]
$(\d C,\d\psi)$ also appears as the boundaries of
\begin{itemize}
\item
an $n$--dimensional $2\epsilon'$ Poincar\'e $2\epsilon'$ quadratic structure
$\Psi_1=(0,\d\psi)$ on the $\epsilon'$ chain map $i_C\co \d C \to C^{n-*}$, and
\item
an $n$--dimensional $\epsilon'$ quadratic structure $\Psi_2=(0,\d\psi)$
on the $0$ chain map $0\co \d C\to 0$, which is $\epsilon'$ Poincar\'e because
$\d C$ is $4\epsilon'$ contractible.
\end{itemize}

The union $\Psi_2\cup_{\d C}-\Psi_3$ is a $600\delta''$ Poincar\'e $7\delta''$
quadratic structure on $0\cup_{\d C}\d D=\C(\d f)$.  By \fullref{boundary}(2), there is a
$2400\delta''$ chain equivalence between $\d D$ and an $n$--dimensional
$500\delta''$ chain complex $\hat\d D$.  This chain equivalence, together with
the $4\epsilon'$ chain contraction of $\d C$, induces a $43200\delta''$ chain
equivalence $g\co 0\cup_{\d C}\d D\to \hat\d D$.  Define a $43200\delta''$
Poincar\'e $3\cdot43200\delta''$ quadratic structure $\hat\psi$ on $\hat\d D$
by $g_\%(\Psi_2\cup_{\d\psi}-\Psi_3)$.  By \fullref{i-str}, there is a $43200\delta''$
Poincar\'e $3\cdot43200\delta''$ quadratic structure on a $43200\delta''$
chain map
\[
(0\cup_{\d C}\d D)\oplus \hat\d D \lga \hat\d D,
\]
and, therefore, the right square in the picture below is filled with a cobordism.

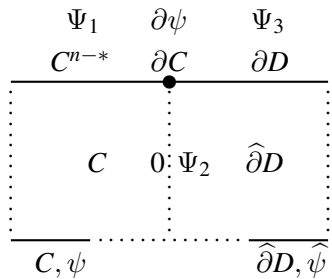
\begin{center}
\begin{picture}(120, 110)(0,-20)
\thicklines
\put(0,0){\line(1,0){30}}
\multiput(34,0)(4,0){14}{\circle*{1}}
\put(90,0){\line(1,0){30}}
\put(0,60){\line(1,0){120}}
\put(60,60){\circle*{5}}
\multiput(0,4)(0,4){14}{\circle*{1}}
\multiput(60,4)(0,4){14}{\circle*{1}}
\multiput(120,4)(0,4){14}{\circle*{1}}
\put(8,-12){$C,\psi$}
\put(92,-12){$\hat\d D,\hat \psi$}
\put(28,26){$C$}
\put(52,26){$0$}
\put(62,26){$\Psi_2$}
\put(88,26){$\hat\d D$}
\put(14, 64){$C^{n-*}$}
\put(52, 64){$\d C$}
\put(90, 64){$\d D$}
\put(20, 80){$\Psi_1$}
\put(52, 80){$\d\psi$}
\put(90, 80){$\Psi_3$}
\end{picture}
\end{center}

The left square can also be filled in with a cobordism.
There is a $3\epsilon'$ homotopy equivalence:
\[
(C^{n-*}\cup_{\d C}0=\C(i_C),\Psi_1\cup_{\d\psi}-\Psi_2) \lga (C,\psi),
\]
and again by \fullref{i-str}, this induces a $30\epsilon'$ Poincar\'e $9\epsilon'$
quadratic structure on a $3\epsilon'$ chain map
\[
(C^{n-*}\cup_{\d C}0) \oplus C \lga C.
\]
Glue these along the pair $(\d C \to 0, \Psi_2)$, and we get a chain map
\[
(C^{n-*}\cup_{\d C} \d D) \oplus C \oplus \hat\d D \lga C \oplus \hat\d D
\]
and a $43200000\delta''$ Poincar\'e $6\cdot 43200\delta''$ quadratic structure
on it.  Since $\d C$ is $4\epsilon'$ contractible and $\d D$ is $2400\delta''$
chain equivalent to $\hat\d D$, there is a $43200\delta''$ chain equivalence
\[
G\co C^{n-*}\cup_{\d C}\d D \lga E=C^{n-*}\oplus \hat\d D,
\]
and hence, by \fullref{i-str3}, there is a $30\cdot43300000\delta''$ Poincar\'e
$4\cdot43300000\delta''$ null-cobordism of $(E,G_\%(\Psi_1\cup_{\d\psi}-\Psi_3))
\oplus(C,-\psi)\oplus(\hat\d D,-\hat\psi)$.
Therefore
\[
[C,\psi]+[\hat\d D,\hat\psi]=[E,G_\%(\Psi_1\cup_{\d\psi}-\Psi_3)]
\]
in $L_n^{13\cdot10^8\delta''}(X;p_X,R)$.

On the other hand, there is a $600\delta''$ Poincar\'e null-cobordism of
$\Psi_1\cup_{\d\psi}-\Psi_3$ on the chain map
\[
C^{n-*}\cup_{\d C}\d D \lga \C(f)^{n+1-*}.
\]
Using $G$ and \fullref{i-str3}, we obtain a $30(600+43200+4)\delta''$
Poincar\'e null-cobordism
\[
(E\to \C(f)^{n+1-*},\quad (\chi,G_\%(\Psi_1\cup_{\d\psi}-\Psi_3)),
\]
and this implies
\[
[E,G_\%(\Psi_1\cup_{\d\psi}-\Psi_3)]= 0 \in L_n^{13\cdot10^8\delta''}(X;p_X,R)
\]
and hence
\[
[C,\psi] = -[\hat\d D,\hat\psi] \in L_n^{13\cdot10^8\delta''}(X;p_X,R).
\]
Since $\d D$ is $400\epsilon''$ contractible over $X-W^{100\epsilon''}$ and
$\hat\d D$ is $2400\delta''$ chain equivalent to $\d D$, $\hat\d D$ is
$5200\delta''$ contractible over $X-W^{2500\delta}$, by \fullref{l-cont}.
By \fullref{5.1-5.2}, there is a $(6n+15)\cdot5200\delta''$ chain equivalence
$h$ from $(\hat\d D,1)$ to an $n$--dimensional $(3n+12)\cdot5200\delta''$
projective chain complex $(F,p)$ on $p_{W^{(20800n+75300)\delta''}}$.
Suppose $\lambda_n \ge 10^5(4n+50)$.
If $\delta\ge \lambda_n\delta''$ and $Y\supset W^{\lambda_n\delta''}$,
then $((F,p),h_\%(\hat\psi))$ represents an element of $L_n^{\{X\},\delta}(Y;p_Y,R)$
by \fullref{i-str}, and its image
\[
i_*([(F,p),h_\%(\hat\psi)])\in L_n^{\alpha\delta}(X;p_X,R)
\]
is represented by $(\hat\d D, (h^{-1})_\%(h_\%(\hat\psi))= (h^{-1}h)_\%(\hat\psi))$.
Since $h^{-1}h$ is $2\delta$ chain homotopic to the identity map,
\[
[\hat\d D,\hat\psi] = [\hat\d D,(h^{-1}h)_\%(\hat\psi)] \in L_n^{\alpha\delta}(X;p_X,R).
\]
Since $\alpha\delta\ge 13\cdot 10^8\delta''$, we have
\[
i_*(-[(F,p),g_\%(\hat\psi)])=[C,\psi] \in L_n^{\alpha\delta}(X;p_X,R).
\]

(3)\qua Let $[C,\psi]\in L_n^{\delta,\epsilon}(X;p_X,R)$, then
$\d C$ is $4\epsilon$ contractible.  Thus $(\d C,1)$ is $4\epsilon$
chain equivalent to $(E=0, q=0)$.

(4)\qua Let $[C,\psi]\in L_n^{\delta',\epsilon'}(X,Y;p_X,R)$ such
that $\d[C,\psi]=0$ in
$L_{n-1}^{\{X\},\delta'',\epsilon''}(W;p_X,R)$.
Let $\beta \co (\d C,1)\to (E,q)$ be a $(200n+300)\epsilon'$ chain equivalence to
an $(n{-}1)$--di\-men\-sion\-al $(100n+300)\epsilon'$ projective chain complex
on $p_{W}$ posited in the definition of $\d$.
By assumption, $[(E,q),\beta_\%\d\psi]=0$ in
$L_{n-1}^{\{X\},\delta'',\epsilon''}(W;p_X,R)$.
By \fullref{g-str2}, there is an $(n{-}1)$--dimensional $100\delta''$ Poincar\'e $2\delta''$
null-cobordism on $p_{W}$
\[
(f'\co (E,q)\lga(D,p),\quad (\delta\psi',\beta_\%\d\psi))
\]
such that $[D,p]=0$ in $\wtilde K_0^{n,9\delta''}(X;p_X,R)$.
By \fullref{l-cont}${}'$, $(\delta\psi',\d\psi)$ is a $125\delta''$ Poincar\'e
$2\delta''$ quadratic structure on the $3\delta''$ chain map
\[
f=f'\circ\beta\co (\d C,1) \lga (D,p).
\]
On the other hand, $(0,\d\psi)$ is a $2\epsilon'$ Poincar\'e $2\epsilon'$
quadratic structure on the $\epsilon'$ chain map
\[
i_C\co (\d C,1) \lga (C^{n-*},1).
\]
Gluing these together, we obtain a $12500\delta''$ Poincar\'e $4\delta''$
quadratic structure
\[
\psi'=(0,\d\psi)\cup_{\d\psi}-(\delta\psi',\d\psi)
\]
on $(C',p')=(C^{n-*},1)\cup_{(\d C,1)}(D,p)$.
Since $n\ge 1$, $(D,p)$ is $540\delta''$ chain equivalent to
an $n$--dimensional free $270\delta''$ chain complex $(F,1)$ on $p_X$
by \fullref{zero-in-K}.

Assume $n\ge2$.
In this case $\d C$ is $12\epsilon'$ chain equivalent to an
$(n{-}1)$--dimensional $4\epsilon'$ chain complex $\hat\d C$, by \fullref{boundary}.
Using these chain equivalences and \fullref{cont-triad}, we can construct a
$6528\delta''$ chain equivalence
\[
\gamma\co (C',p') \lga (C''=C^{n-*}\cup_{\hat\d C} F,1).
\]
If $\delta\ge 9\cdot 10^5\delta''$, then
$(C'', \psi''=\gamma_\%\psi')$ determines an element of $L_n^\delta(X;p_X,R)$.
Suppose $\delta\ge 4\cdot 10^6\delta''$ and $Y'\supset W^{12\cdot 10^6\delta''}$.
We shall show that its image by $j_*$ is equal to the relax-control image of
$[C,\psi]$ in $L_n^\delta(X,Y';p_X,R)$.

Since $(D,p)$ lies over $W$, it is $0$ contractible over $X-W$.
Therefore, by \fullref{cont-triad}, the chain map $G\co (C',p')\to (C^{n-*}\cup_{\d C}0,1)$
defined by
\[
G=\begin{pmatrix}1 & 0 & 0\\ 0 & 1 & 0\end{pmatrix}\co (C^{n-r},1)\oplus(\d C_{r-1},1)\oplus
(D_r,p)\lga (C^{n-r},1)\oplus(\d C_{r-1},1)
\]
is a $18\delta''$ chain equivalence over $X-W^{6\delta''}$.
Furthermore one can easily check that $G$ is 0 connected and that
$G_\%(\psi')=(0,\d\psi)\cup_{\d \psi}-(0,\d\psi)$.
Compose $G$ with a $3\epsilon'$ homotopy equivalence
\[
((C^{n-*}\cup_{\d C}0,1),(0,\d\psi)\cup_{\d \psi}-(0,\d\psi))\lga((C,1),\psi)
\]
to get a $3\epsilon'$ connected $19\delta''$ chain equivalence over
$X-W^{7\delta''}$:
\[
H\co (C',p')\lga (C,1);\qquad H_\%(\psi')=\psi.
\]
By \fullref{i-str}, there is a $3\cdot 12500\delta''$ connected $3\cdot 19\delta''$
quadratic structure $(0,\psi'\oplus -\psi)$ on a chain map
$(C',p')\oplus (C,1)\to(C,1)$ that is $125000\delta''$ Poincar\'e over
$X-W^{375007\delta''}$.
Now use the $6528\delta''$ chain equivalence $\gamma\co (C',p')\to (C'',1)$ and
\fullref{i-str3}${}'$ to this cobordism to obtain an $(n{+}1)$--dimensional $\delta$ cobordism
between $(C'',\psi'')$ and $(C,\psi)$ that is $\delta$ connected and
$\delta$ Poincar\'e over $X-Y'$.

In the $n=1$ case, use the non-positive chain complex
obtained from $\d C$ by applying the folding argument from top
instead of $\hat\d C$.
See the proof of \fullref{boundary}(1).

(5)\qua Let $[C,\psi]\in L_n^{\delta,\epsilon}(X,Y;p_X,R)$ and let $\beta\co (\d C,1)\to (E,q)$
be as in the definition of $\d$; $\d[C,\psi]$ is given by
$[(E,q),\beta_\%\d\psi]$.
There exist $(n{-}1)$--dimensional free $\epsilon'$ chain complexes
$E'$, $F$ on $p_X$ and a $3\epsilon'$ chain equivalence
\[
g\co (E,q)\oplus(\d\Sigma E',1) \lga (F,1)
\]
with $\Sigma E'$ given the trivial quadratic structure $\theta$, and
$i_*[(E,q),\beta_\%\d\psi]$ is represented by
$(F,g_\%((\beta_\%\d\psi)\oplus \d\theta))$.
We construct a $\delta''$ Poincar\'e null-cobordism of this.

Take the direct sum of the algebraic Poincar\'e thickenings of $(C,\psi)$
and $(\Sigma E',\theta)$ to get an $\epsilon'$ Poincar\'e pair
\[
(\d C\oplus \d\Sigma E'\lga C^{n-*}\oplus E'^{n-1-*},\quad (0\oplus0,\d\psi
\oplus \d\theta)).
\]
Apply the $4\epsilon'$ chain equivalence
\[
\d C\oplus \d\Sigma E'=(\d C,1)\oplus(\d\Sigma E',1)\xrightarrow{\beta\oplus1}{}
(E,q)\oplus(\d\Sigma E',1)\xrightarrow{~g~}{}(F,1)=F,
\]
to this pair, to obtain an $\epsilon''$ Poincar\'e null-cobordism
of $(F,g_\%((\beta_\%\d\psi)\oplus \d\theta)))$.
(If $n\ge 2$, then we may assume $E'=0$, and the proof can be much simplified.)

(6)\qua Take an element $[(C,p),\psi]\in L_{n-1}^{\{X\},\delta',\epsilon'}(Y;p_X,R)$ and assume
$i_*[(C,p),\psi]=0$ in $L_{n-1}^{\delta'',\epsilon''}(X;p_X,R)$.
By definition of $i_*$, there exist $(n{-}1)$--dimensional free $\epsilon'$
chain complexes $E$, $F$ on $p_X$ and a $3\epsilon'$ chain equivalence
$g\co (C,p)\oplus(\d\Sigma E,1)\to (F,1)$ such that $i_*[(C,p),\psi]=
[F,g_\%(\psi\oplus\d\theta)]$.  Here $\theta$ is the trivial quadratic
structure on $\Sigma E$.  By \fullref{g-str}, there is an $n$--dimensional
$100\delta''$ Poincar\'e $2\delta''$ null-cobordism on $p_X$ of
$(F,g_\%(\psi\oplus\d\theta))$:
\[
(f\co F\lga D, \quad (\delta\psi, g_\%(\psi\oplus\d\theta))).
\]
By \fullref{l-cont}${}'$, we obtain a $127\delta''$ Poincar\'e $3\delta''$
null-cobordism:
\[
(f\circ g\co (C,p)\oplus(\d\Sigma E,1)\lga(D,1),\quad (\delta\psi,\psi\oplus\d\theta)).
\]
Take the union of this with the 0 connected $\epsilon'$ projective quadratic pair
\[
((C,p)\lga 0,\quad (0,-\psi)),
\]
which is 0 Poincar\'e over $X-Y$, and the $3\epsilon'$ Poincar\'e $3\epsilon'$
quadratic pair
\[
(i_{\Sigma E}\co (\d\Sigma E,1)\lga (E^{n-*},1),\quad (0,-\d\theta))
\]
to get a $6\delta''$ projective quadratic complex $((\hat C,\hat r),\hat\psi)$
which is $12700\delta''$ Poincar\'e over $X-Y^{12700\delta''}$ and is $12700\delta''$
connected.

The $3\epsilon'$ chain equivalence $g$ induces a $48\delta''$ chain equivalence
$\tilde g\co (\hat C,\hat r)\to (\tilde C,1)$ to an $n$--dimensional free chain
complex $(\tilde C,1)=(D,1)\cup_{(F,1)}(E^{n-*},1)$,
and the $18\delta''$ quadratic structure $\tilde\psi=\tilde g_\%(\hat\psi)$ is
$2\cdot10^5\delta''$ Poincar\'e over $X-Y^{4\cdot10^5\delta''}$ and is
$2\cdot10^5\delta''$ connected.  Suppose $W\supset Y^{10^6\delta''}$ and
$\delta\ge10^6\delta''$.  Then $(\tilde C,\tilde\psi)$ defines an
element in $L_n^\delta(X,W;p_X,R)$.

We shall show that $\d[\tilde C,\tilde\psi]=[(C,p),\psi]$ in
$L_{n-1}^{\{X\},\kappa_n\delta}(W^{\kappa_n\delta};p_X,R)$.
By the definition of $\d$, there is a $(200n+300)\delta$ chain equivalence
$\beta\co (\d\tilde C,1)\to(\tilde E,\tilde q)$ to an $(n{-}1)$--dimensional
$(100n+300)\delta$ projective chain complex on $p_{W^{(150n+300)\delta}}$,
and $\d[\tilde C,\tilde\psi]$ is represented by
$((\tilde E,\tilde q),\beta_\%\d\tilde\psi)$.
We construct a cobordism between $((C,p),\psi)$
and $((\tilde E,\tilde q),\beta_\%\d\tilde\psi)$.

By \fullref{i-str}${}'$, $\tilde g$ induces an $(n{+}1)$--dimensional
$3\cdot48\delta''$ cobordism:
\[
(\begin{pmatrix}\tilde g & 1\end{pmatrix}:
(\hat C,\hat r)\oplus (\tilde C,1)\lga(\tilde C,1),
\quad \tilde\Psi=(0,\hat\psi\oplus-\tilde\psi)).
\]
Let us apply the boundary construction to this to get a $6\cdot48\delta''$
chain map
\[
\d\begin{pmatrix}\tilde g&1\end{pmatrix}:
\d(\hat C,\hat r)\oplus \d(\tilde C,1)\lga (G,q)
\]
and a $9\cdot48\delta''$ Poincar\'e $6\cdot48\delta''$ quadratic structure
$(\chi, \d\hat\psi\oplus-\d\tilde\psi)$ on it.  We modify this to get
the desired cobordism.

Firstly, note that $((\hat C,\hat r),\hat\psi)$ is the algebraic Thom complex
of a $12700\delta''$ Poincar\'e $6\delta''$ quadratic pair with boundary equal
to $((C,p),\psi)$.  Therefore there is a $11\cdot12700\delta''$ chain equivalence
$\gamma\co \d(\hat C,\hat r)\to (C,p)$ such that $\gamma_\%(\d\hat\psi)=\psi$.

Secondly, there is a chain equivalence
$\beta\co (\d\tilde C,1)\to (\tilde E, \tilde q)$ as noted above.

Thirdly, recall that $(G,q)$ is equal $\Omega\C(\D_{\tilde\Psi})$ and
$\d(\tilde C,1)$ is equal to $\Omega\C(\D_{\tilde\psi})$, and note that
there is a $96\delta''$ chain equivalence
\[
\C(\begin{pmatrix}\tilde g & 1\end{pmatrix})^{n+1-*}
\xrightarrow{(0~~1~~-\tilde g^*)}{} (\hat C,\hat r)
\xrightarrow{({\tilde g}^{-1})^*}{} (\tilde C,1)
\]
and that it induces a $6\delta$ chain equivalence from $(G,q)$ to
$\d(\tilde C,1)$.  Compose this with $\beta$ to get a $(200n+306)\delta$
chain equivalence $\beta'\co (G,q)\to (\tilde E,\tilde q)$.

Now, by \fullref{i-str3}${}'$, one can conclude that the chain equivalences
$\gamma$, $\beta$, $\beta'$ induce an $n$--dimensional $\kappa_n\delta$ Poincar\'e
cobordism on $p_{W^{\kappa_n\delta}}$:
\[
((C,p)\oplus(\tilde E,\tilde q)\lga(\tilde E,\tilde q),\quad (\chi,\psi\oplus
-\beta_\%(\d\tilde\psi))).
\]
Since $[C,p]=0$ in $\wtilde K_0^{n,\epsilon'}(X;p_X,R)$ and
$[\tilde E,\tilde q]=[\d\tilde C,1]=0$ in
$\wtilde K_0^{n,\kappa_n\delta}(X;p_X,R)$,
this implies that $[(C,p),\psi]=\d[\tilde C,\tilde\psi]$ in
$L_{n-1}^{\{X\},\kappa_n\delta}(W^{\kappa_n\delta};p_X,R)$.
\end{proof}

\section{Excision}\label{Exc}

\showlabel{Exc}
In this section we study the excision property of epsilon-controlled $L$--theory.
Suppose that $X$ is the union of two closed subsets
$A$ and $B$ with intersection $M=A\cap B$.
There is an inclusion-induced homomorphism
\[
i_*\co L_n^{\delta,\epsilon}(A,M;p_{A},R) \to L_n^{\delta,\epsilon}(X,B;p_X,R).
\]
For $n\ge 1$, we construct its stable inverse
\[
\exc\co L_n^{\delta,\epsilon}(X,B;p_X,R)\to
L_n^{\delta,\epsilon}(A,A\cap M^{(n+5)4\delta};p_{A},R) .
\]
First we define geometric subcomplexes and quotient complexes of
free chain complexes.  Let $C$ be a free chain complex on $p_X$.
When each $C_r$ is the direct sum $C_r=C'_r\oplus
C''_r$ of two geometric submodules and $d_C$ is of the form
\[
\begin{pmatrix}d_{C'} & * \cr 0 & d_{C''}\end{pmatrix}
\co C'_r\oplus C''_r \lga C'_{r-1}\oplus C''_{r-1},
\]
$C'$ is said to be a \emph{geometric subcomplex} of $C$,
and $C''$ (together with $d_{C''}$) is said to be the \emph{quotient} of $C$
by $C'$ and is denoted $C/C'$.
If $C$ is a free $\epsilon$ chain complex, then any geometric subcomplex
$C'$ and the quotient $C/C'$ are both free $\epsilon$ chain complexes.
The obvious projection map $p\co C \to C/C'$ is 0 connected.

Next suppose we are given an $n$--dimensional $\epsilon$ quadratic
complex $(C,\psi)$ on $p_X$ and $C'$ is a geometric subcomplex of $C$.
The projection $p\co C\to C/C'$ induces an $n$--dimensional $\epsilon$
quadratic complex $(C/C',p_\%\psi)$ and there is an $\epsilon$ cobordism
between $(C,\psi)$ and $(C/C',p_\%\psi)$.
For a morphism $g\co G\to H$ between geometric modules and
geometric submodules $G'\subset G$ and $H'\subset H$,
we write $g(G')\subset H'$ when every path with non-zero coefficient in $g$
starting from a generator of $G'$ ends at a generator of $H'$.

\begin{proposition}\label{quotients}
\showlabel{quotients}
Let $(C,\psi)$, $C'$, and $p$ be as above, and
suppose $(C,\psi)$ is $\epsilon$ connected.
If $\D_\psi(C'^n)\subset C'_0$, then $(C/C',p_\%\psi)$ and the cobordism
between $(C,\psi)$ and $(C/C',p_\%\psi)$ induced by $p$ are both
$\epsilon$ connected.
\end{proposition}

\begin{proof}
Let us write $C''=C/C'$.
By assumption, the morphism $d_{\C(\D_\psi)}\co \C(\D_\psi)_1 \to \C(\D_\psi)_0$
can be expressed by a matrix of the form
\[
\begin{pmatrix}d_{C'} & * & * & *\\ 0 & d_{C''} & 0 & \D_{p_\%\psi}\end{pmatrix}
\co C'_1\oplus C''_1 \oplus C'^n \oplus C''^n \to C'_0 \oplus C''_0.
\]
Let $h\co C_0=\C(\D_\psi)_0\to \C(\D_\psi)_1$ be a $4\epsilon$ morphism
such that $d_{\C(\D_\psi)}h\sim_{8\epsilon}1_{C_0}$, and define
$4\epsilon$ morphisms $h_1\co C''_0\to C''_1$
and $h_2\co C''_0 \to C''^n$ by
\[
h=\begin{pmatrix}* & *\\ * & h_1\\ * & *\\ * & h_2\end{pmatrix}\co 
C'_0\oplus C''_0 \to C'_1\oplus C''_1 \oplus C'^n \oplus C''^n.
\]
Then we get a homotopy
\[
d_{C''}h_1 + \D_{p_\%\psi}h_2~ \sim_{8\epsilon}~ 1_{C''_0}.
\]
Therefore
\[
\begin{pmatrix}h_1\cr h_2\end{pmatrix}\co C''_0 \sa C''_1\oplus C''^n
\]
gives a desired splitting of the boundary morphism
$\C(\D_{p_\%\psi})_1 \to \C(\D_{p_\%\psi})_0$.
Therefore $(C'', p_\%\psi)$ is $\epsilon$ connected.
Now the $\epsilon$ connectivity of cobordism induced by $p$
follows from \fullref{conn}.
\end{proof}

\begin{example}\label{ex-subcx}
\showlabel{ex-subcx}
Let $(C,\psi)$ be an $n$--dimensional $\epsilon$ quadratic
complex on $p_X$ and $Y$ be a subset of $X$.  Fix $\delta(>0)$ and $l(\ge 0)$,
and define a geometric submodule $C'_r$ of $C_r$ to be the
restriction
$C_r(Y^{(n+l-r)\delta})$ of $C_r$ to $Y^{(n+l-r)\delta}$.
If $\delta\ge \epsilon$, $\{C'_r\}$ is a geometric subcomplex of $C$,
and we can form the quotient $C/C'$ of $C$ by $C'$ and the natural
projection $p\co C\to C/C'$.
$(C/C')_r$ is equal to $C_r(X-Y^{(n+l-r)\delta})$.
Suppose further that $(C,\psi)$ is $\epsilon$ connected,
$\delta\ge4\epsilon$, and $n\ge 1$; then $\D_\psi(C'^n)\subset C'_0$
holds, and $(C/C',p_\%\psi)$ and the cobordism between $(C,\psi)$
and $(C/C',p_\%\psi)$ induced by $p$ are both $\epsilon$ connected.
\end{example}

Next we consider pairs.
Suppose $(f\co C\to D,(\delta\psi,\psi))$ is an $(n{+}1)$--dimensional
$\epsilon$ quadratic pair on $p_X$ and $C'$, $D'$ are geometric subcomplexes
of $C$, $D$, respectively such that $f(C'_r)\subset (D'_r)$ for every $r$.
Define an $\epsilon$ chain map $\wbar f\co C/C'\to D/D'$ by
\[
f=\begin{pmatrix}* & * \\ 0 & \wbar f\end{pmatrix}\co 
C'_r\oplus (C/C')_r \sa D'_r\oplus (D/D')_r,
\]
then the diagram
\[
\xymatrix{C \ar[d]_p \ar[r]^f & D \ar[d]^q \\
C/C' \ar[r]_{\wbar f} & D/D'}
\]
commutes strictly, where $p$ and $q$ are the natural projections, and
\[
(\wbar f\co C/C'\to D/D',(q_\%\delta\psi, p_\%\psi))
\]
is an $(n{+}1)$--dimensional $\epsilon$ quadratic pair.

\begin{proposition}\label{inducedmap}
\showlabel{inducedmap}
If $(f,(\delta\psi,\psi))$ is $\epsilon$ connected, $\D_{\psi}(C'^n)\subset C'_0$,
and $\D_{\delta\psi}(D'^{n+1})\subset D'_0$, then
$(\wbar f, (q_\%\delta\psi, p_\%\psi))$ is $\epsilon$ connected.
\end{proposition}

\begin{proof}
We check the $\epsilon$ connectivity of the duality map
$\D_{(q_\%\delta\psi, p_\%\psi)}$.
Let us use the notation $C''=C/C'$ and $D''=D/D'$.
The boundary morphism
$$d_{\C(\D_{(\delta\psi,\psi)})}\co 
\C(\D_{(\delta\psi,\psi)})_1 \to \C(\D_{(\delta\psi,\psi)})_0$$
can be expressed by a matrix of the form
\begin{multline*}
\begin{pmatrix}* & * & * & * & * & * \\
  0 & d_{D''} & 0 & \D_{q_\%\delta\psi} & 0 & \wbar f\D_{p_\%\psi}
\end{pmatrix}\co  \\
D'_1\oplus D''_1\oplus D'^{n+1}\oplus D''^{n+1}\oplus C'^n\oplus C''^n
\sa ~D'_0\oplus D''_0.
\end{multline*}
The desired $\epsilon$ connectivity follows from this as in
\fullref{quotients}.
\end{proof}

\begin{proposition}\label{compare}
\showlabel{compare}
Let $Y$ be a subset of $X$, and let
$[(C,d),\psi]$ and $[(\hat C,\hat d),\hat\psi]$ be elements
of $L_n^{\delta,\epsilon}(X,Y;p_X,R)$ $(n\ge 1)$.
If
\begin{enumerate}
\item $C_r(X-Y)=\hat C_r(X-Y)$,
\item $d_r|X-Y^{4\epsilon}=\hat d_r|X-Y^{4\epsilon}$, and
\item $\psi_s|X-Y^{4\epsilon}=\hat \psi_s|X-Y^{4\epsilon}$
\end{enumerate}
for every $r$ and $s~(\ge0)$, then $[(C,d),\psi]=[(\hat C,\hat d),\hat\psi]$
in $L_n^{\delta,\epsilon}(X,Y^{(n+3)4\epsilon};p_X,R)$.
\end{proposition}

\begin{proof}
Define  a geometric subcomplex $C'$ of $C$ by
$C'_r=C_r(Y^{(n+1-r)4\epsilon})$, and let $p\co C\to C/C'$ be the projection.
Then $(C/C',p_\%\psi)$ is an $\epsilon$ connected $\epsilon$ quadratic
complex by \fullref{quotients}.  The boundary maps for $\C(\D_\psi)$ have radius $4\epsilon$
and are of the form
\begin{multline*}
\begin{pmatrix}
d_{C'} & * & * & *\\
0 & d_{C''} & * & (-)^{r-1}\D_{p_\%\psi}\\
0 & 0 & * & * \\
0 & 0 & * & (-)^{r-1}d_{C''}^*
\end{pmatrix}
\co C'_r \oplus C''_r \oplus C'^{n+1-r} \oplus C''^{n+1-r} \\
\sa C'_{r-1}\oplus C''_{r-1} \oplus C'^{n+2-r} \oplus C''^{n+2-r}.
\end{multline*}
Therefore $\C(\D_{p_\%\psi})$ and $\C(\D_\psi)$ are exactly the same
over $X-Y^{(n+2)4\epsilon}$, and
$\C(\D_{p_\%\psi})$ is $4\epsilon$ contractible over $X-Y^{(n+3)4\epsilon}$,
ie $p_\%\psi$ is $\epsilon$ Poincar\'e over $X-Y^{(n+3)4\epsilon}$.
In fact, if $\Gamma$ is a $4\epsilon$ chain contraction over $X-Y$ of
$\C(\D_\psi)$, then $\Gamma|X-Y^{(n+2)4\epsilon}$ gives a $4\epsilon$
chain contraction over $X-Y^{(n+3)4\epsilon}$ of $\C(\D_{p_\%\psi})$.
Thus $(C/C',p_\%\psi)$ determines an element of
$L_n^{\delta,\epsilon}(X,Y^{(n+3)4\epsilon};p_X,R)$.

By \fullref{conn}, the cobordism between $(C,\psi)$ and $(C/C',p_\%\psi)$ induced
by $p$ is an $\epsilon$ connected $\epsilon$ quadratic pair.
Since this cobordism is exactly the same over $X-Y^{(n+2)4\epsilon}$ as
the trivial cobordism from $(C,\psi)$ to itself, it is $\epsilon$
Poincar\'e over $X-Y^{(n+3)4\epsilon}$.  Therefore,
\[
[C,\psi]~=~[C/C',p_\%\psi]~\in ~L_n^{\delta,\epsilon}(X,Y^{(n+3)4\epsilon};p_X,R).
\]
The same construction for $(\hat C,\hat\psi)$ yields the same
element as this, and we can conclude that
\[
[C,\psi]~=~[\hat
C,\hat\psi]~\in~L_n^{\delta,\epsilon}(X,Y^{(n+3)4\epsilon};p_X,R).\proved
\]
\end{proof}

Now suppose $X$ is the union of two
closed subsets $A$, $B$ with intersection $N=A\cap B$.

\begin{lemma}\label{effect}
\showlabel{effect}
Let $G$, $H$ be geometric modules on $p_X$, and
$f\co G\to H$ be a morphism of radius $\delta$.  Then, for any $\gamma\ge0$,
\[
f(G(B\cup N^\gamma))~\subset~
  H(B\cup N^{\max\{\gamma+\delta,2\delta\}}).
\]
\end{lemma}

\begin{proof}
This can be deduced from the following two claims:
\begin{enumerate}
\item $f(G(N^\gamma))~\subset~ H(N^{\gamma+\delta})$,
\item $f(G(B))~\subset~ H(B\cup N^{2\delta})$.
\end{enumerate}
The first claim is obvious.  To prove the second claim,
take a generator of $G(B)$  and a path $c$ starting from $a$
with non-zero coefficient in $f$.
By its continuity, the path $p_X\circ c$ in $X$ either stays inside of
$B$ or passes through a point in $N$, and hence its image is
contained in $B\cup N^{2\delta}$.  This proves the second claim.
\end{proof}

Now let us define the excision map:
\[
\exc\co L_n^{\delta,\epsilon}(X,B;p_X,R)\to
L_n^{\delta,\epsilon}(A,A\cap N^{(n+5)4\delta};p_{A},R),
\]
Take an element $[C,\psi]\in L_n^{\delta,\epsilon}(X,B;p_X,R)$.
Define a geometric subcomplex $C'$ of $C$ by
\[
C'_r~=~C_r(B\cup N^{(n+2-r)4\epsilon}),
\]
and let $p\co C\to C/C'$ denote the projection.
Then, by \fullref{effect} and \fullref{quotients}, $(C/C',p_\%\psi)$ is an
$\epsilon$ connected $\epsilon$ quadratic complex on $p_{A}$ and is $\epsilon$
Poincar\'e over $A-N^{(n+4)4\epsilon}$.
We define $\exc([C,\psi])$ to be the element
\[
[C/C',p_\%\psi]~\in ~L_n^{\delta,\epsilon}(A,A\cap N^{(n+5)4\epsilon};p_{A},R).
\]
The excision map is well-defined. Suppose
\[
[C,\psi]~=~[\hat C,\hat\psi]~\in~ L_n^{\delta,\epsilon}(X,B;p_X,R).
\]
Without loss of generality we may assume that there is a $\delta$
connected $\delta$ cobordism
\[
(f\co C\oplus \hat C \to D, (\delta\psi, \psi\oplus -\hat\psi))
\]
between $(C,\psi)$ and $(\hat C,\hat\psi)$
that is $\delta$ Poincar\'e over $X-B$.
Let us now construct $(C/C',p_\%\psi)$ and $(\hat C/\hat C', \hat p_\%\hat\psi)$
as above, define a geometric subcomplex $D'$ of $D$ by
\[
D'_r~=~D_r(B\cup N^{(n+3-r)4\delta}),
\]
and let $q\co D\to D/D'$ denote the projection.
By \fullref{effect} and \fullref{inducedmap}, we obtain an $\delta$ connected $\delta$
cobordism
\[
(\wbar f \co C/C'\oplus \hat C/{\hat C}' \to D/D', (q_\%\delta\psi,
p_\%\psi\oplus -\hat p_\%\hat\psi))
\]
which is $\delta$ Poincar\'e over $A-B\cup N^{(n+5)4\delta}$.
Therefore $\exc$ is well-defined.

By using \fullref{compare},
we can check that the homomorphisms $i_*$ and $\exc$ are stable inverses;
ie the following diagram commutes:
\[
\xymatrix{L_n^{\delta,\epsilon}(A,N;p_{A},R)\ar[d] \ar[r]^{i_*} &
L_n^{\delta,\epsilon}(X,B;p_X,R) \ar @{=} [d]\\
L_n^{\delta,\epsilon}(A,A\cap N^{(n+5)4\delta};p_{A},R)\ar@{=}[d]
& \ar[l]_(.4){\rm exc} L_n^{\delta,\epsilon}(X,B;p_X,R) \ar[d]\\
L_n^{\delta,\epsilon}(A,A\cap N^{(n+5)4\delta};p_{A},R) \ar[r]_-{i_*} &
L_n^{\delta,\epsilon}(X,B\cup N^{(n+5)4\delta};p_X,R)}
\]
\noindent
where the vertical maps are the homomorphisms induced by inclusion maps.

\section{Mayer--Vietoris type sequence}\label{MV}

\showlabel{MV}
We continue to assume that  $X$ is the union of two
closed subsets $A$, $B$ with intersection $N=A\cap B$,
and will present a Mayer--Vietoris type stably exact sequence.

Replace $\kappa_n$ by $\kappa_n +4(n+5)$, and
suppose $\delta\ge\epsilon>0$.
Let $W$ be a subset of $X$ containing $N^{\kappa_n\delta}$
and assume $\delta'\ge\kappa_n\delta$, $\epsilon'\ge\kappa_n\epsilon$.
Then a homomorphism
\[
\wbar\d\co L_n^{\delta,\epsilon}(X;p_X,R)\sa
L_{n-1}^{\{A\cup W\},\delta',\epsilon'}(W;p_{A\cup W},R)
\]
is obtained by composing the following maps:
\begin{multline*}
L_n^{\delta,\epsilon}(X;p_X,R)~ \sa ~L_n^{\delta,\epsilon}(X,B;p_X,R)~
\xrightarrow{\rm{exc}}~
L_n^{\delta,\epsilon}(A,A\cap N^{(n+5)4\delta};p_{A},R)\\
\xrightarrow{~\d~} ~
L_{n-1}^{\{A\},\delta',\epsilon'}(A\cap W;p_{A},R) ~\sa ~
L_{n-1}^{\{A\cup W\},\delta',\epsilon'}(W;p_{A\cup W},R)
\end{multline*}
If $[C,\psi]\in L_n^{\delta,\epsilon}(X;p_X,R)$, then its image $\wbar\d[C,\psi]$
is represented by $((E,q),\psi')$ which is homotopy equivalent to
the boundary $\hat\d(C/C',p_\%\psi)$, where $C'\subset C$ and $p\co C\to C/C'$ are
as in the definition of \exc.
This is exactly the projective quadratic Poincar\'e complex $(Q,\wbar\psi)$ which appears
in the Splitting Lemma:

\begin{lemma}[Pedersen--Yamasaki \cite{P-Y}]
\label{splitlemma}
\showlabel{splitlemma}
For any integer $n\ge 2$,
there exists a positive number $\kappa_n\ge 1$ such that the following holds:
Suppose $p_X\co M\to X$ is a map to a metric space $X$,
$X$ is the union of two closed subsets $A$ and $B$ with intersection $N=A\cap B$,
and $R$ is a ring with involution.  Let $\epsilon$ be any positive number, and
set $\epsilon'=\kappa_n\epsilon$, $N'=N^{\epsilon'}$, $A'=A\cup N'$, and
$B'=B\cup N'$.  Then for any $n$--dimensional quadratic Poincar\'e $R$--module complex
$c=(C,\psi)$ on $p_X$ of radius $\epsilon$,
there exist a Poincar\'e cobordism of radius $\epsilon'$
from $c$ to the union
$c'\cup c''$ of an $n$--dimensional quadratic Poincar\'e pair
$c'=(f'\co Q\to C', (\delta\wbar\psi', -\wbar\psi))$ on $p_{A'}$ of radius $\epsilon'$ and
an $n$--dimensional quadratic Poincar\'e pair
$c''=(f''\co Q\to C'', (\delta\wbar\psi'', \wbar\psi))$
on $p_{B'}$ of radius $\epsilon'$ along an $(n{-}1)$--dimensional
quadratic Poincar\'e projective $R$--module complex $(Q,\wbar\psi)$ on $p_{N'}$,
where $Q$ is $\epsilon'$ chain equivalent to an
$(n{-}1)$--dimensional free chain complex on $p_{A'}$
and also to an
$(n{-}1)$--dimensional free chain complex on $p_{B'}$.
\end{lemma}

From this and its relative version, we obtain the following:

\begin{proposition}\label{MVboundary}
\showlabel{MVboundary}
If $n\ge 2$, the map $\wbar\d$ factors through
a homomorphism
\[
\d\co L_n^{\delta,\epsilon}(X;p_X,R)\sa L_{n-1}^{{\F},\delta',\epsilon'}(W;p_X,R),
\]
where ${\F} = \{A\cup W,A\cup W\}$.
Moreover the image $\d[C,\psi]$ is given by $[Q,\wbar\psi]$ which appears
in any splitting (up to cobordism) of $(C,\psi)$ according to the closed subsets $A$, $B$.
\end{proposition}

Now we present the Mayer--Vietoris type stably-exact sequence.
It is made up of three kinds of maps.
The first is the map
\[
i_*\co L_n^{\{A,B\},\delta,\epsilon}(N;p_X,R)
\to L_n^{\delta',\epsilon'}(A;p_A,R)\oplus L_n^{\delta',\epsilon'}(B;p_B,R)
\]
defined by $i_*(x)=(\iota_A(x), -\iota_B(x))$
when $\delta'\ge\alpha\delta$ and $\epsilon'\ge\alpha\epsilon$.
The second is the map
\[
j_*\co 
L_n^{\delta,\epsilon}(A;p_A,R)\oplus L_n^{\delta,\epsilon}(B;p_B,R)
\to L_n^{\delta',\epsilon'}(X;p_X,R)
\]
defined by $j_*(x,y)=j_{A*}(x)+j_{B*}(y)$
when $\delta'\ge\delta$ and $\epsilon'\ge\epsilon$.
Here $j_A\co A\to X$ and $j_B\co B\to X$ are inclusion maps.
The third is the map $\d$ given in \fullref{MVboundary}:
\[
\d\co L_n^{\delta,\epsilon}(X;p_X,R)\sa
L_{n-1}^{\{A\cup W,B\cup W\},\delta',\epsilon'}(W;p_X,R),
\]
where $W\supset N^{\kappa_n\delta}$,
$\delta'\ge\kappa_n\delta$, and $\epsilon'\ge\kappa_n\epsilon$.

In the rest of this section, we omit the control map and the coefficient ring
from notation.

\begin{theorem}\label{ThMVS1}
\showlabel{ThMVS1}
For any integer $n\ge 2$, there exists a constant $\lambda_n > 1$
which depends only on $n$ such that the following holds true for
any control map $p_X$ and two closed subsets $A$, $B$ of $X$ satisfying
$X=A\cup B$:

{\rm (1)}\qua Suppose $\delta'\ge \alpha\delta$, $\epsilon'\ge \alpha\epsilon$,
$\delta''\ge\delta'$, and $\epsilon''\ge\epsilon'$ so that
the following two maps are defined:
\[
L_n^{\{A,B\},\delta,\epsilon}(N)\xrightarrow{~i_*~}
L_n^{\delta',\epsilon'}(A)\oplus L_n^{\delta',\epsilon'}(B)
\xrightarrow{~j_*~}
L_n^{\delta'',\epsilon''}(X)
\]
Then $j_*i_*$ is zero.

{\rm (2)}\qua
Suppose $\delta''\ge\delta'$, $\epsilon''\ge\epsilon'$ so that
$j_*\co L_n^{\delta',\epsilon'}(A)\oplus L_n^{\delta',\epsilon'}(B)
\to L_n^{\delta'',\epsilon''}(X)$
is defined.
If $\delta\ge\lambda_n\delta''$ and $W\supset N^{\lambda_n\delta''}$,
then the relax-control image of the kernel of $j_*$ in
$L_n^{\alpha\delta}(A\cup W)\oplus L_n^{\alpha\delta}(B\cup W)$ is
contained in the image of $i_*$ below:
\[
\xymatrix{
& L_n^{\delta',\epsilon'}(A)\oplus L_n^{\delta',\epsilon'}(B)
\ar[r]^-{j_*~~} \ar[d]
& L_n^{\delta'',\epsilon''}(X) \\
L_n^{\{A\cup W, B\cup W\},\delta}(W) \ar[r]^-{i_*}
& L_n^{\alpha\delta}(A\cup W)\oplus L_n^{\alpha\delta}(B\cup W)
& \\
}
\]
{\rm (3)}\qua
Suppose $\delta'\ge\delta$, $\epsilon'\ge\epsilon$,
$W\supset N^{\kappa_n\delta'}$, $\delta''\ge\kappa_n\delta'$,
and $\epsilon''\ge\kappa_n\epsilon'$
so that the following two maps are defined:
\[
L_n^{\delta,\epsilon}(A)\oplus L_n^{\delta,\epsilon}(B)
\xrightarrow{j_*}{}
L_n^{\delta',\epsilon'}(X)
\xrightarrow{\ \d\ }{}
L_{n-1}^{\{A\cup W,B\cup W\},\delta'',\epsilon''}(W)
\]
Then $\d j_*$ is zero.

{\rm (4)}\qua
Suppose $W\supset N^{\kappa_n\delta'}$, $\delta''\ge \kappa_n\delta'$,
and $\epsilon''\ge\kappa_n\epsilon'$ so that the map
$\d\co L_n^{\delta',\epsilon'}(X) \to
L_{n-1}^{\{A\cup W,B\cup W\},\delta'',\epsilon''}(W)$
is defined.
If $\delta\ge\lambda_n\delta''$, then the relax-control
image of the kernel of $\d$ in $L_n^{\delta}(X)$ is
contained in the image of $j_*$ below:
\[
\xymatrix{
 & L_n^{\delta',\epsilon'}(X)
\ar[r]^-{\d} \ar[d]
 & L_{n-1}^{\{A\cup W,B\cup W\},\delta'',\epsilon''}(W)
\\
L_n^{\delta}(A\cup W)\oplus L_n^{\delta}(B\cup W)
\ar[r]^-{j_*}
 & L_n^{\delta}(X)
& \\
}
\]

{\rm (5)}\qua Suppose $W\supset N^{\kappa_n\delta}$, $\delta'\ge\kappa_n\delta$,
$\epsilon'\ge\kappa_n\epsilon$,
$\delta''\ge\alpha\delta'$, and $\epsilon''\ge\alpha\epsilon'$
so that the following two maps are defined:
\[
L_n^{\delta,\epsilon}(X)
\xrightarrow{\ \d\ }
L_{n-1}^{\{A\cup W,B\cup W\},\delta',\epsilon'}(W)
\xrightarrow{i_*}
L_{n-1}^{\delta'',\epsilon''}(A\cup W)\oplus L_{n-1}^{\delta'',\epsilon''}(B\cup W)
\]
Then $i_*\d$ is zero.

{\rm (6)}\qua
Suppose $\delta''\ge \alpha\delta'$, and $\epsilon''\ge\alpha\epsilon'$ so that
$i_*\co 
L_{n-1}^{\{A,B\},\delta',\epsilon'}(N)\to
L_{n-1}^{\delta'',\epsilon''}(A)\oplus L_{n-1}^{\delta'',\epsilon''}(B)$
is defined.
If $\delta\ge\lambda_n\delta''$, $N'\supset N^{\lambda_n\delta''}$,
and $W=(N')^{\kappa_n\delta}$, then the relax-control
image of the kernel of $i_*$ in $L_{n-1}^{\{A\cup W,B\cup W\},\delta}(W)$ is
contained in the image of $\d$ associated with the two closed subsets
$A\cup N'$, $B\cup N'$:
\[
\xymatrix{
 & L_{n-1}^{\{A,B\},\delta',\epsilon'}(N)
\ar[r]^-{i_*} \ar[d]
 & L_{n-1}^{\delta'',\epsilon''}(A)\oplus L_{n-1}^{\delta'',\epsilon''}(B)
\\
L_n^{\delta}(X)
\ar[r]^-{\d}
 & L_{n-1}^{\{A\cup W,B\cup W\},\delta}(W)
& \\
}
\]
\end{theorem}

\begin{proof}
(1)\qua Take an element $x=[Q,\psi]\in L_n^{\{A,B\},\delta,\epsilon}(N)$.
The image $i_*(x)$ is a pair $([c_A], -[c_B])$ where $c_A$ and $c_B$ are
free quadratic Poincar\'e complexes on $p_A$ and $p_B$ that are both
homotopy equivalent to $(Q,\psi)$, and hence $[c_A]=[c_B]\in L_n^{\delta'',\epsilon''}(X)$.
Therefore, $j_*i_*(x)=[c_A]-[c_B] = 0$.

(2)\qua First, temporarily use the constant
$\lambda_n$ posited in the splitting lemma. Take an element
$x=([C_A,\psi_A],[C_B,\psi_B])\in L_n^{\delta',\epsilon'}(A)\oplus
L_n^{\delta',\epsilon'}(B)$ and assume $j_*(x)=0$. There exists a
null-cobordism $(f\co C_A\oplus C_B \to D, (\delta\psi,\psi_A\oplus
-\psi_B))$. Its boundary is already split according to $A$ and
$B$, so use the relative splitting to this null-cobordism to get
cobordisms of radius $\lambda_n\delta''$:
\begin{gather*}
(f_A\co (C_A,1)\oplus Q\to (D_A,1),(\delta\psi_A,\psi_A\oplus -\wbar\psi))
~~\hbox{on}~~p_{A\cup W},\\
(f_B\co (C_B,1)\oplus Q\to (D_B,1),(\delta\psi_B,\psi_B\oplus -\wbar\psi))
~~\hbox{on}~~p_{B\cup W}.
\end{gather*}
By the Poincar\'e duality $D_A^{n+1-*}\simeq \C(f_A)$, we have
\[
[D_A]-[C_A]-[Q]=[\C(f_A)]=[D_A^{n+1-*}]=0\in\wtilde K_0^{n+1,4\lambda_n\delta''}(A\cup W),
\]
and, hence, we have $[Q]=0$ in $\wtilde K_0^{n,36\lambda_n\delta''}(A\cup W)$.
See Ranicki--Yamasaki \cite[Section~3]{RY-K}.
Similarly $[Q]=0$ in $\wtilde K_0^{n,36\lambda_n\delta''}(B\cup W)$.
Thus we obtain an element
$[Q,\wbar\psi]$ of $L_n^{\{A\cup W, B\cup W\},36\lambda_n\delta''}(W)$.
Replace $\lambda_n$ by something bigger (at least $36\lambda_n$) so that
its image via $i_*$ in $L_n^{\alpha\delta}(A\cup W)\oplus L_n^{\alpha\delta}(B\cup W)$
is equal to $([C_A,\psi_A],[C_B,\psi_B])$ whenever $\delta\ge\lambda_n\delta''$.

(3)\qua If we start with an element $x=([C_A,\psi_A], [C_B,\psi_B])$, then
$j_*(x)$ is represented by $(C_A,\psi_A)\oplus(C_B,\psi_B)$ which is already
split according to $A$ and $B$.  Therefore $\d j_*(x)=0$.

(4)\qua Temporarily set the constant $\lambda_n$
to be the one posited in the splitting lemma. Take an element
$[C,\psi]$ in $L_n^{\delta',\epsilon'}(X)$ such that
$\d[C,\psi]=0$. $(C,\psi)$ splits into two adjacent pairs:
\[
a=(f_A\co Q\to (C_A,1), (\delta\psi_A, -\wbar\psi))~~~\hbox{and}~~~
b=(f_B\co Q\to (C_B,1), (\delta\psi_B, \wbar\psi))
\]
such that $[Q,\wbar\psi]=0$ in $L_{n-1}^{\{A\cup W,B\cup W\},\delta'',\epsilon''}(W)$.
Take a $\delta''$ null-cobordism on $p_W$ $p=(g\co Q\to P, (\delta\wbar\psi,\wbar\psi))$
such that the reduced projective class of $P$ is zero on $p_{A\cup W}$ and also
on $p_{B\cup W}$.
$C_A\cup_Q P$ is chain equivalent to an $n$--dimensional free complex $F_A$ on $p_A$, and
$C_B\cup_Q P$ is chain equivalent to an $n$--dimensional free complex $F_B$ on $p_B$.
Use these to fill in the bottom squares with cobordisms:

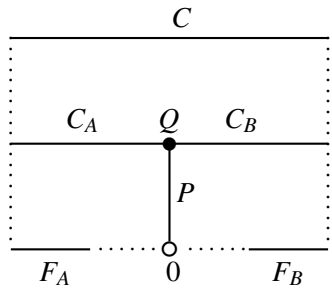
\begin{center}
\begin{picture}(120, 110)(0,-20)
\thicklines
\put(0,0){\line(1,0){30}}
\multiput(34,0)(4,0){6}{\circle*{1}}
\multiput(68,0)(4,0){6}{\circle*{1}}
\put(90,0){\line(1,0){30}}
\put(0,40){\line(1,0){120}}
\put(0,80){\line(1,0){120}}
\put(60,40){\circle*{5}}
\put(60,0){\circle{5}}
\multiput(0,4)(0,4){20}{\circle*{1}}
\multiput(120,4)(0,4){20}{\circle*{1}}
\put(60,3){\line(0,1){35}}
\put(10,-12){$F_A$}
\put(58,-12){$0$}
\put(98,-12){$F_B$}
\put(62,18){$P$}
\put(20, 46){$C_A$}
\put(55, 46){$Q$}
\put(80, 46){$C_B$}
\put(60, 84){$C$}
\end{picture}
\end{center}

Replacing $\lambda_n$ with something larger if necessary, we obtain
free quadratic Poincar\'e complexes on $p_{A\cup W}$ and $p_{B\cup W}$
whose sum is $\lambda_n\delta''$ cobordant to $(C,\psi)$.

(5)\qua If we start with an element $x=[C,\psi]$ in $L_n^{\delta,\epsilon}(X)$,
then $j_*(x)$ is represented by the projective piece $(Q,\wbar\psi)$ obtained
by splitting, and the null-cobordisms required to show $i_*\d(x)=0$ are easily
constructed from the split pieces.

(6)\qua
Take an element $[Q,\psi]$ of $L_{n-1}^{\{A,B\},\delta',\epsilon'}(N)$.
Then $(Q,\psi)$ is homotopy equivalent to a free quadratic Poincar\'e complex
$((C_A,1), \psi_A)$ on $p_A$ and also to a free quadratic Poincar\'e complex
$((C_B,1), \psi_B)$ on $p_B$.
If $i_*[Q,\psi]=0$, then these are both null-cobordant;
there are quadratic Poincar\'e pairs
\begin{align*}
&(f_A\co C_A\to D_A, (\delta\psi_A, \psi_A)) \hbox{~~~on $p_A$, and}\\
&(f_B\co C_B\to D_B, (\delta\psi_B, \psi_B)) \hbox{~~~on $p_B$.}
\end{align*}
Use the homotopy equivalence $(C_A,\psi_A)\simeq (C_B,\psi_B)$
to replace the boundary of the latter by $(C_A,\psi_B)$,
and glue them together to get an element $[D,\delta\psi]$ of $L_n^{\delta}(X)$
for some $\delta>0$.
Note that $(D,\delta\psi)$ has a splitting into two pairs with the common boundary
piece equal to $(Q,\psi)$, so we have $\d[D,\delta\psi]=[Q,\psi]$ in
$L_{n-1}^{\{A\cup W,B\cup W\},\delta}(W)$.
\end{proof}

\section{A special case}\label{SpC}

\showlabel{SpC}
In this section we treat the case when there are no controlled $K$--theoretic
difficulties.

First assume that $X$ is a finite polyhedron.
We fix its triangulation.
Under this assumption we can simplify the Mayer--Vietoris type sequence
of the previous section at least for sufficiently small $\epsilon$'s and $\delta$'s.
$X$ is equipped with a deformation $\{f_t\co X\to X\}$
called `rectification' \cite{P-Y} which
deforms sufficiently small neighborhoods of the $i$--skeleton $X^{(i)}$ into $X^{(i)}$
such that $f_t$'s are uniformly Lipschitz.
This can be used to rectify the enlargement of the relevant subsets at the expense of
enlargement of $\epsilon$'s and $\delta$'s.
We thank Frank Quinn for showing us his description of {\it uniformly continuous CW complexes}
which are designed for taking care of these situations in a more general setting.

Next let us assume that $X$ is a finite polyhedron and that
the control map $p_X\co M\to X$ is a fibration with path-connected fiber $F$ such that
\[
\Wh(\pi_1(F)\times {\ZZ}^k) = 0
\]
for every $k\ge 0$.
The condition on the fundamental group is satisfied if $\pi_1(F)\cong \ZZ^l$ for some $l\ge 0$.
If we study the proofs of \cite[8.1 and 8.2]{RY-K} carefully, we obtain the following.

\begin{proposition}\label{CKvanish}
\showlabel{CKvanish}
Let $p_X$ be as above and  $n\ge 0$ be an integer.
Then there exist numbers $\epsilon_0>0$ and $0<\mu\le 1$ which depend on $X$ and $n$
such that the relax-control maps
\begin{align*}
&\wtilde K_0^{n,\epsilon}(S;p_S,\ZZ) \lga \wtilde K_0^{n,\epsilon'}(S;p_S,\ZZ)\\
&\Wh^{n,\epsilon}(S;p_S,\ZZ) \lga \Wh^{n,\epsilon'}(S;p_S,\ZZ)
\end{align*}
are zero maps for any subpolyhedron $S$, any $\epsilon'\le\epsilon_0$
and any $\epsilon\le\mu\epsilon'$.
\end{proposition}

This means that there is a homomorphism functorial with respect to relaxation of control
\[
L_n^{\F,\delta,\epsilon}(S;p_X,\ZZ)\lga
L_n^{\F\cup\{S\},\delta',\epsilon'}(S;p_X,\ZZ)
\]
for any family $\F$ of subpolyhedra of $X$ containing $S$, if
$\delta'\le\epsilon_0$,  $\delta\le\mu\delta'$,  and  $\epsilon\le\mu\epsilon'$.
Compose this with the homomorphism
\[
\iota_S \co L_n^{\F\cup\{S\},\delta',\epsilon'}(S;p_X,\ZZ) \lga
L_n^{\alpha\delta',\alpha\epsilon'}(S;p_X,\ZZ)
\]
to get a homomorphism
\[
\iota\co L_n^{\F,\delta,\epsilon}(S;p_X,\ZZ)\lga
L_n^{\alpha\delta',\alpha\epsilon'}(S;p_X,\ZZ).
\]
A stable inverse $\tau$ functorial with respect to relaxation of control can be defined
by $\tau([C,\psi])=[(C,1),\psi]$, and we have a commutative diagram:
\[
\xymatrix{L_n^{\delta,\epsilon}(S;p_{X},\ZZ)\ar[d] \ar[r]^{\tau} &
L_n^{\F, \delta,\epsilon}(S;p_X,\ZZ) \ar @{=} [d]\\
L_n^{\alpha\delta',\alpha\epsilon'}(S;p_{X},\ZZ)\ar@{=}[d]
& \ar[l]_(.4){\iota} L_n^{\F,\delta,\epsilon}(S;p_X,\ZZ) \ar[d]\\
L_n^{\alpha\delta',\alpha\epsilon'}(S;p_{X},\ZZ) \ar[r]_-{\tau} &
L_n^{\F,\alpha\delta',\alpha\epsilon'}(S;p_X,\ZZ)}.
\]
Thus the Mayer--Vietoris type sequence is stably exact
when we replace the controlled projective $L$--group terms
with appropriate controlled $L$--groups.

Furthermore, since $p_X$ is a fibration, we have a stability for controlled
$L$--groups:

\begin{proposition}[Pedersen--Yamasaki {\cite[Theorem~1]{P-Y}}]\label{stability}
\showlabel{stability}
Let $n\ge 0$.
Suppose $Y$ is a finite polyhedron and $p_Y\co M\to Y$ is a fibration.
Then there exist constants $\delta_0>0$ and $K>1$,
which depends on the integer $n$ and  $Y$, such that the relax-control map
$L_n^{\delta',\epsilon'}(Y;p_Y,R)\to L_n^{\delta,\epsilon}(Y;p_Y,R)$ is an isomorphism
if $\delta_0\ge \delta\ge K\epsilon$, $\delta_0\ge \delta'\ge K\epsilon'$,
$\delta\ge\delta'$,  and $\epsilon\ge\epsilon'$.
\end{proposition}

Now let us denote these isomorphic groups $L_n^{\delta,\epsilon}(Y;p_Y,R)$
($\delta_0\ge\delta$, $\delta\ge\kappa\epsilon$) by $L_n(Y;p_Y,R)$.
When the coefficient ring $R$ is $\ZZ$, we omit $\ZZ$ and use the notation
$L_n(Y;p_Y)$.

\begin{theorem}\label{periodicity}
\showlabel{periodicity}
Let $p_X\co M\to X$ be a fibration over a finite polyhedron $X$.
Then $L_n(X;p_X,R)$ is 4--periodic:
$L_n(X;p_X,R)\cong L_{n+4}(X;p_X,R)$ ($n\ge 0$).
\end{theorem}

\begin{proof}
The proof of the 4--periodicity of $L_n({\mathbb A})$ of an additive category with involution
given in Ranicki \cite{AdditiveL} adapts well to the controlled setting.
\end{proof}

We have a Mayer--Vietoris exact sequence for $L_n$ with coefficient ring $\ZZ$.

\begin{theorem}\label{ThMVS2}
\showlabel{ThMVS2}
Let $X$ be a finite polyhedron and suppose that $p_X\co M\to X$
is a fibration with path-connected fiber $F$ such that
$\Wh(\pi_1(F)\times {\ZZ}^k) = 0$ for every $k\ge 0$.
If $X$ is the union of two subpolyhedra $A$ and $B$,
then there is a long  exact sequence
\begin{gather*}
\dots \xrightarrow{~\d~} L_n({A\cap B};p_{A\cap B})\xrightarrow{~i_*~}
L_n(A;p_A)\oplus L_n(B;p_B)
\xrightarrow{~j_*~} L_n(X;p_X)~\\
\xrightarrow{~\d~} L_{n-1}({A\cap B};p_{A\cap B})\xrightarrow{~i_*~} \dots
\xrightarrow{~j_*~} L_0(X;p_X).
\end{gather*}
\end{theorem}

\begin{proof}
The exactness at the term $L_2(A;p_A)\oplus L_2(B;p_B)$ and
at the terms to the left of it follows immediately from the stably-exact sequence.
The exactness at other terms follows from the 4--periodicity.
\end{proof}

Recall that there is a functor $\LL(-)$ from spaces to $\Omega$--spectra
such that $\pi_n(\LL(M))=L_n(\ZZ[\pi_1(M)])$ constructed geometrically by Quinn \cite{Q1},
and algebraically by Ranicki \cite{Topman}.
Blockwise application of $\LL$ to $p_X$ produces a generalized homology group
$H_n(X;\LL(p_X))$ (see Quinn \cite{Q2}).
There is a map $A\co H_n(X;\LL(p_X))\to L_n(X;p_X)$ called the assembly map.
See Yamasaki \cite{MY87} for the $\LL^{-\infty}$--analogue, involving the lower $L$--groups
of Ranicki \cite{LowerKL}.

\begin{theorem}\label{assembly}
\showlabel{assembly} Let $X$ be a finite polyhedron and suppose
that $p_X\co M\to X$ is a fibration with path-connected fiber $F$
such that $\Wh(\pi_1(F)\times {\ZZ}^k) = 0$ for every $k\ge 0$.
Then the assembly map $A\co H_n(X;\LL(p_X))\to L_n(X;p_X)$ is an
isomorphism.
\end{theorem}

\begin{proof}
We actually prove the isomorphism $A\co H_n(S;\LL(p_S))\to L_n(S;p_S)$
for all the subpolyhedra $S$ of $X$ by induction on the number of simplices.

When $S$ consists of a single point $v$, then the both sides are
$L_n(\ZZ[\pi_1(p_X^{-1}(v))])$
and $A$ is the identity map.

Suppose $S$ consists of $k>1$ simplices and assume by induction
that the assembly map is an isomorphism for all subpolyhedra consisting of less
number of simplices.
Pick a simplex $\Delta$ which is not a face of other simplices and let $A=\Delta$ and
$B=S-\hbox{interior}(\Delta)$.
Since $A$ contracts to a point $v$, it can be easily shown that
$H_n(A;\LL(p_A))$ and $L_n(A;p_A)$ are both $L_n(\ZZ[\pi_1(p_X^{-1}(v))])$, and
the assembly map $A\co  H_n(A;\LL(p_A))\to L_n(A;p_A)$ is an isomorphism.
By induction hypothesis the assembly maps for $B$ and $A\cap B$ are both isomorphisms.
We can conclude that the assembly map for $S$ is an isomorphism by an application
of 5--lemma to the ladder made up of the Mayer--Vietoris sequences
for $H_*(-)$ and $L_*(-)$.
\end{proof}

\begin{remarkun}\label{assembly2}
\showlabel{assembly2}
If $F$ is simply-connected, then
$\Wh(\pi_1(F)\times {\ZZ}^k) = \Wh({\ZZ}^k)=0$ for every $k\ge 0$
by the celebrated result of Bass, Heller and Swan. In this case
$H_n(X;\LL(p_X))$ is isomorphic to
the generalized homology group $H_n(X;\LL)$ where $\LL$ is the
4--periodic simply-connected surgery spectrum with $\pi_n(\LL)=L_n(\ZZ[\{1\}])$ and
we have an assembly isomorphism
$$A\co H_n(X;\LL)\cong L_n(X;p_X).$$
This is the controlled surgery obstruction group which appears in the
controlled surgery exact sequence of Pedersen--Quinn--Ranicki \cite{PQR} (as required for the
surgery classification of exotic homology manifolds in
Bryant--Ferry--Mio--Weinberger \cite{BFMW}).
There the control map is not assumed to be a fibration.  We believe
that most of the arguments in this paper work in a more general
situation.
\end{remarkun}

As an application of \fullref{ThMVS2}, we consider the $\ZZ$--coefficient controlled $L$--group of
$p_X\times 1\co M\times S^1 \to X\times S^1$.

\begin{corollary}\label{xcircle}
\showlabel{xcircle} Let $n\ge 0$, and
let $X$ and $p_X\co M\to X$ be as in \fullref{ThMVS2}.  Then there is
a split short exact sequence
\[
0\to L_n(X;p_X) \xrightarrow{i_*} L_n(X\times S^1;p_X\times 1)\xrightarrow{B} L_{n-1}(X;p_X)\to 0.
\]
\end{corollary}

\begin{proof}
Split the circle $S^1=\d([-1,1]\times[-1,1])\subset {\RR}^2$ into two pieces
\[
S^1_+=\{(x,y)\in S^1|y\ge 0\} \quad \hbox{ and } \quad S^1_-=\{(x,y)\in S^1|y\le 0\},
\]
with intersection $\{p=(1,0),~q=(-1,0)\}$. Let $\d$ be the
connecting homomorphism in the Mayer--Vietoris sequence
\fullref{ThMVS2} corresponding to this splitting, and consider the
composite
\begin{multline*}
B\co L_n(X\times S^1;p_X\times 1)\xrightarrow{\d}
L_{n-1}(X\times\{p\};p_X\times 1)\oplus L_{n-1}(X\times\{q\};p_X\times 1)\\
\xrightarrow{\hbox{\scriptsize projection}}
L_{n-1}(X\times\{p\};p_X\times 1)\cong L_{n-1}(X;p_X).
\end{multline*}
Then $\d$ can be identified with
\[
(B, -B)\co L_n(X\times S^1;p_X\times 1)\lga L_{n-1}(X;p_X)\oplus L_{n-1}(X;p_X).
\]
The map $i_*$ is the map induced by the inclusion map
\[
L_n(X;p_X)\cong L_n(X\times\{p\};p_X\times 1)\xrightarrow{i_*}
L_n(X\times S^1;p_X\times 1).
\]
The exactness follows easily from the exactness of the Mayer--Vietoris sequence.
A splitting of $B$ can be constructed by gluing two product cobordisms.
\end{proof}

\begin{corollary}\label{torus}
\showlabel{torus}
Let $T^n$ be the $n$--dimensional torus $S^1\times \dots \times S^1$.
Then
\begin{align*}
L_m(T^n;1_{T^n})&\cong
\bigoplus_{r=0}^n\begin{pmatrix}n\\r\end{pmatrix}L_{m-r}(\ZZ)\\
&\cong L_m(\ZZ[\pi_1(T^n)])\qquad\hbox{\rm ($m\ge n$)}.
\end{align*}
\end{corollary}

\begin{proof}
Use \fullref{xcircle} repeatedly to obtain
$$L_m(T^n;1_{T^n})\cong
\bigoplus_{r=0}^n\begin{pmatrix}n\\r\end{pmatrix}L_{m-r}(\ZZ).$$
The isomorphism
$$\bigoplus_{r=0}^n\begin{pmatrix}n\\r\end{pmatrix}L_{m-r}(\ZZ)\cong
L_m(\ZZ[\pi_1(T^n)])$$
is the well-known computation obtained geometrically by Shaneson and Wall,
and algebraically by Novikov and Ranicki.
\end{proof}


\section{Locally finite analogues}\label{LFA}
\showlabel{LFA}
Up to this point, we considered only finitely generated modules and chain complexes.
In this section we deal with infinitely generated objects;
such objects arise naturally when we take the pullback of
a finitely generated object via an infinite-sheeted covering map.
We restrict ourselves to a very special case necessary for our application.

\begin{definition}[{{Ranicki and Yamasaki \cite[page 14]{RY-K}}}]
Consider the product $M\times N$ of two spaces.  A geometric module on
$M\times N$ is said to be $M$--\emph{finite} if, for any $y\in N$, there
is a neighbourhood $U$ of $y$ in $N$ such that $M\times U$ contains only
finitely many basis elements; a projective module $(A,p)$ on $M\times N$
is said to be $q$--\emph{finite} if $A$ is $M$--finite; a projective
chain complex $(C,p)$ on $M\times N$ is $M$--\emph{finite} if each
$(C_r,p_r)$ is $M$--finite.  (In \cite{RY-K}, we used the terminology
``$M$--\emph{locally finite}'', but this does not sound right and we
decided to use ``$M$\emph{--finite}'' instead.)  When $M$ is compact,
$M$--finiteness is equivalent to the ordinary locally-finiteness.
\end{definition}

Consider a control map $p_X\co M\to X$ to a metric space $X$, and
let $N$ be another metric space.
Give the maximum metric to the product $X\times N$, and let us use the map
\[
p_X \times 1_N\co M\times N \sa X\times N,
\]
as the control map for $M\times N$.

\begin{definition}
For $\delta\ge\epsilon>0$, $Y\subset X$, and a family $\F$ of subsets of $X$ containing $Y$,
define {\it $M$--finite $(\delta,\epsilon)$--controlled $L$--groups}
$L_n^{M,\delta,\epsilon}(X\times N, Y\times N;p_X\times 1,R)$,
and {\it $M$--finite $(\delta,\epsilon)$--controlled projective $L$--groups}
$L_n^{M,\F,\delta,\epsilon}(Y\times N; p_X\times 1,R)$
by requiring that every chain complexes concerned are $M$--finite.
\end{definition}

All the materials up to \fullref{MV} are valid for $M$--finite analogues.
In the previous section, there are several places where we assumed $X$ to be a finite polyhedron,
and they may not automatically generalize to the $M$--finite case.

The most striking phenomenon about $M$--finite objects is the following vanishing result
on the half line.

\begin{proposition}\label{vanish}
\showlabel{vanish}
Let $p_X\co M\to X$ be a control map, $N$ a metric space,
and give $N\times [0,\infty)$ the maximum metric.
For any $\epsilon>0$  and $\delta\ge \epsilon$,
\begin{align*}
&L_n^{M,\delta,\epsilon}(X\times N \times [0,\infty); p_X\times 1,R) = 0, \\
&\wtilde K_0^{M,n,\epsilon}(X\times N\times [0,\infty); p_X\times 1, R) = 0.
\end{align*}
\end{proposition}

\begin{proof}
This is done using repeated shifts towards infinity and the `Eilenberg Swindle'.
Let us consider the case of
$L_n^{M,\delta,\epsilon}(X\times N\times [0,\infty);p_X\times 1,R)$.
Let $J=[0,\infty)$ and define $T\co M\times N\times J\to M\times N \times J$
by $T(x,x',t)=(x,x',t+\epsilon)$.
Take an element $[c]\in L_n^{M,\delta,\epsilon}(X\times N\times J,p_X\times 1,R)$.
It is zero, because there exist $M$--finite $\epsilon$ Poincar\'e cobordisms:
\begin{align*}
c ~&\sim~ c\oplus (T_\#(-c) \oplus T_\#^2(c)) \oplus (T_\#^3(-c) \oplus
T_\#^4(c)) \oplus \cdots\\
&=~ (c\oplus T_\#(-c)) \oplus (T_\#^2(c)\oplus T_\#^3(-c)) \oplus \cdots
~\sim~ 0.
\end{align*}
The proof for controlled $\wtilde K$ is similar.  See the appendix to \cite{RY-K}.
\end{proof}

Thus, the analogue of Mayer--Vietoris type sequence \eqref{ThMVS1}
for the control map $p_X\times 1\co M\times N\times \RR\to X\times N\times \RR$
with respect to the splitting
$X\times N\times \RR=X\times N\times (-\infty,0]\cup X\times N\times [0,\infty)$
reduces to
\[
0~\sa~L_n^{M,\delta,\epsilon}(X\times N\times \RR;p_X\times 1,R)\xrightarrow{~\d~}
L_{n-1}^{M,p,\delta',\epsilon'}(X\times N\times I; p_X\times 1,R)~\sa~0,
\]
where $\delta'=\kappa_n\delta$, $\epsilon'=\kappa_n\epsilon$,
$I = [-\kappa_n\delta,\kappa_n\delta]$, and
the right hand side is the $M$--finite projective $L$--group
$L_{n-1}^{M,\{~\},\delta',\epsilon'}(X\times N\times I; p_X\times 1,R)$
corresponding to the empty family $\F=\{~\}$.

A diagram chase shows that there exists a well-defined homomorphism
\[
\beta\co L_{n-1}^{M,p,\delta',\epsilon'}(X\times N\times I; p_X\times 1,R)~\sa~
L_n^{M,\delta'',\epsilon''}(X\times N\times \RR;p_X\times 1,R),
\]
where $\gamma''=\lambda_n\kappa_n\lambda_{n-1}\alpha\gamma'$ and
$\epsilon''=\lambda_n\kappa_n\lambda_{n-1}\alpha\epsilon'$.
The homomorphisms $\d$ and $\beta$ are stable inverses of each other;
the compositions are both relax-control maps.

Note that, for any $\delta\ge\epsilon>0$, the retraction induces an isomorphism
\[
L_{n-1}^{M,p,\delta,\epsilon}(X\times N\times I; p_X\times 1,R)~\cong~
L_{n-1}^{M,p,\delta,\epsilon}(X\times N\times \{0\}; p_X,R).
\]
Thus, we have obtained:

\begin{theorem}\label{split1}
\showlabel{split1}
Splitting along $X\times N\times \{0\}$  induces a stable isomorphism:
\[
\d\co L_n^{M,\delta,\epsilon}(X\times N\times \RR; p_X\times 1,R)~\lga~
L_{n-1}^{M,p,\delta',\epsilon'}(X\times N; p_X\times 1 ,R).
\]
\end{theorem}

Now , as in the previous section,
let us assume that $X$ is a finite polyhedron and $p_X\co M\to X$ is a fibration
with a path-connected fiber $F$ such that $\Wh(\pi_1(F)\times \ZZ^k)=0$ for every $k\ge 0$.

The following is an $M$--finite analogue of \fullref{CKvanish}.

\begin{proposition}\label{LFCKvanish}
\showlabel{LFCKvanish}
Let $p_X$ be as above and  $n\ge 0$, $k\ge 0$ be integers.
Then there exist numbers $\epsilon_0>0$ and $0<\mu\le 1$ which depend on $X$, $n$, and $k$
such that the relax-control maps
\begin{align*}
&\wtilde K_0^{M,n,\epsilon}(X\times\RR^k;p_X\times 1,\ZZ) \lga
\wtilde K_0^{M,n,\epsilon'}(X\times\RR^k;p_X\times 1,\ZZ)\\
&\Wh^{M,n,\epsilon}(X\times\RR^{k};p_X\times 1,\ZZ) \lga
\Wh^{M,n,\epsilon'}(X\times\RR^{k};p_X\times 1,\ZZ)
\end{align*}
is the zero map for any $\epsilon'\le\epsilon_0$
and any $\epsilon\le\mu\epsilon'$.
\end{proposition}

\begin{proof}
First note that,
since $X\times\RR^k$ is not a finite polyhedron unless $k=0$, the proof for
\fullref{CKvanish} does not directly apply to the current situation.

Let us consider the Whitehead group case first.
Since the $k=0$ case was already treated in \fullref{CKvanish}, let us suppose $k>0$.
Let $T^k$ denote the $k$--torus $(S^1)^k$, and
define $p^{(k)}_X\co M\times T^k\to X$ to be the following composite map:
\[
M\times T^k\xrightarrow{\hbox{\small projection}} M \xrightarrow{p_X} X.
\]
By the Mayer--Vietoris type sequence for controlled $K$--theory, the group
\[
\Wh^{M,n,\epsilon}(X\times\RR^{k};p_X\times 1,\ZZ)
\]
is stably isomorphic to
\[
\wtilde K_0^{M,n-1,\epsilon}(X\times\RR^{k-1};p_X\times 1,\ZZ),
\]
which is also stably isomorphic to
\[
\Wh^{M\times S^1,n,\epsilon}(X\times\RR^{k-1};p'_X\times 1,\ZZ).
\]
The last statement is a locally-finite analogue of \cite[7.1]{RY-K}.
The proof given there works equally well here.
Therefore
$\Wh^{M,n,\epsilon}(X\times\RR^{k};p_X\times 1,\ZZ)$
is stably isomorphic to
\[
\Wh^{M\times T^k,n,\epsilon}(X;p^{(k)}_X,\ZZ),
\]
for which the stable vanishing is already known.
This completes the Whitehead group case.

The $\wtilde K_0$ case follows from the stable vanishing of
\[
\Wh^{M,n+1,\epsilon}(X\times\RR^{k+1};p_X\times 1,\ZZ).
\proved
\]
\end{proof}

From this we get:

\begin{proposition}\label{split2}
\showlabel{split2}
Assume that $X$ is a finite polyhedron and $p_X\co M\to X$ is a fibration
with a path-connected fiber $F$ such that $\Wh(\pi_1(F)\times \ZZ^k)=0$
for every $k\ge 0$.  Splitting along $X\times\RR^{m-1}\times\{0\}$ induces a stable isomorphism
\[
\d\co L_n^{M,\delta,\epsilon}(X\times \RR^{m}; p_X\times 1,\ZZ)~\lga~
L_{n-1}^{M,\delta',\epsilon'}(X\times \RR^{m-1}; p_X\times 1 ,\ZZ).
\]
\end{proposition}

\begin{corollary}\label{LFstability}
\showlabel{LFstability}
Let $X$ and $p_X$ be as above, then stability holds for
$L_n^{M,\delta,\epsilon}(X\times \RR^{m}; p_X\times 1,\ZZ)$;
that is, it is isomorphic to the limit
\[
L_n^{M}(X\times \RR^{m}; p_X\times 1,\ZZ)=
\lim_{0<\epsilon\ll\delta\to 0}L_n^{M,\delta,\epsilon}(X\times \RR^{m}; p_X\times 1,\ZZ)
\]
when $0<\epsilon\ll\delta$ and $\delta$ is sufficiently small.
\end{corollary}

\begin{proof}  By the 4--periodicity, we may assume that $n>m$.
Then the proposition above gives a stable isomorphism with
$L_{n-m}^{\delta,\epsilon}(X; p_X,\ZZ)$, and the result follows.
\end{proof}

\begin{corollary}\label{split3}
\showlabel{split3}
Let $X$ and $p_X$ be as above, then
splitting along $X\times\RR^{m-1}\times\{0\}$ induces an isomorphism
\[
\d\co L_n^{M}(X\times \RR^{m}; p_X\times 1,\ZZ) \lga
L_n^{M}(X\times \RR^{m-1}\times\{0\}; p_X\times 1,\ZZ).
\]
\end{corollary}

\begin{proof}
Immediate from \fullref{split2} and \fullref{LFstability}.
\end{proof}

\section{Controlled surgery obstructions}\label{CSO}
\showlabel{CSO}

\def\lf{\hbox{\it\scriptsize lf}}
We discuss the controlled surgery obstructions and an application.
We only consider the identity control maps on polyhedra or on the products of
polyhedra and $\RR^m$.
$X$--finiteness on $X\times\RR^m$ is the same as the usual local finiteness,
so we use the following notation throughout this section:
\begin{eqnarray*}
L_n^{\lf,\delta,\epsilon}(X\times\RR^m) & = & L_n^{X,\delta,\epsilon}(X\times\RR^m;1_X\times 1,\ZZ)\\
L_n^{\lf}(X\times\RR^m) & = & L_n^{X}(X\times\RR^m;1_X\times 1,\ZZ)
\end{eqnarray*}
We omit the decoration `{\it lf\/}' when $m=0$.

Let $(f,b)\co M\to N$ be a degree 1 normal map between connected
oriented closed $PL$ manifolds of dimension $n$. Quadratic
construction on this produces an element $\sigma_N(f,b)\in
L_n^{\delta,\epsilon}(N)$ for arbitrarily small
$\delta\gg\epsilon>0$ (see Ranicki--Yamasaki \cite{RY-SQC}). By
\fullref{stability}, this
defines an element $\sigma_N(f,b)\in L_n^{\lf}(N)$. This is the
{\it controlled surgery obstruction} of $(f,b)$, and its image via
the forget-control map
\[
L_n(N)\to L_n(\{\hbox{\it pt.}\};N\to \{\hbox{\it pt.}\},\ZZ)=L_n(\ZZ[\pi_1(N)])
\]
is the ordinary surgery obstruction $\sigma(f,b)$ of $(f,b)$.
The controlled surgery obstruction $\sigma_N(f,b)$ vanishes, if $(f,b)$ is
normally bordant to a sufficiently small homotopy equivalence measured on $N$.

Similarly, if $(f,b)\co V\to W$ is a degree 1 normal map between connected open oriented $PL$
manifolds of dimension $n$, we obtain its {\it controlled surgery obstruction}
$\sigma_W(f,b)$ in $L_n^{\lf}(W)$.

\begin{theorem}\label{signature}
\showlabel{signature}
Let $X$ be a connected oriented closed $PL$ manifold of dimension $4k$,
and $f\co V^n\to W^n=X\times\RR^{n-4k}$ be a homeomorphism of open $PL$ manifolds.
Homotope $f$ to produce a map $g\co V\to W$ which is transverse regular to $X\times\{0\}
\subset X\times\RR^{n-4k}$.  Then the $PL$ submanifold $g^{-1}(X\times\{0\})$ of $V$
and $X$ have the same signature:
$\sigma(g^{-1}(X\times\{0\})) = \sigma(X)$.
\end{theorem}

\begin{proof}
The homeomorphism $f$ determines a degree 1 normal map $F$ between $V$ and $W$, and hence
determines an element $\sigma_W(F)\in L_n^{\lf}(X\times\RR^{n-4k})$.
Repeated application of splitting \fullref{split3} induces an isomorphism
\[
\d^{n-4k}\co L_n^{\lf}(X\times\RR^{n-4k})\to L_{4k}(X).
\]
The image of $\sigma_W(F)$ by this map is the controlled surgery obstruction
$\sigma_X(g|,b)$ of the degree 1 normal $PL$ map $(g|,b)\co Y=g^{-1}(X\times\{0\})\to X\times\{0\}=X$.
Since $f$ is a homeomorphism, $F$ is normally bordant to arbitrarily small
homotopy equivalences.  Therefore, $\sigma_W(F)$ is zero and hence $\sigma_X(g|,b)$ is zero.
This means that the ordinary surgery obstruction $\sigma(g|,b)$ is also zero.
The equality $\sigma(Y)=\sigma(X)$ follows from this.
\end{proof}

Now we apply the above to prove the topological invariance of the rational
Pontrjagin classes (see Novikov \cite{Novikov}).

\begin{theorem}[S\,P Novikov]
If $h\co M^n\to N^n$ is a homeomorphism between oriented closed $PL$ manifolds,
then $h^*p_i(N)=p_i(M)$,
where $p_i$ are the rational Pontrjagin classes.
\end{theorem}

\begin{proof}
Recall that the rational Pontrjagin classes $p_*(N)\in H^{4*}(N;\QQ)$ determine and are
determined by the ${\mathcal L}$--genus ${\mathcal L}_*(N)\in H^{4*}(N;\QQ)$,
and that the degree $4k$ component ${\mathcal L}_k(N)\in H^{4k}(N;\QQ)$ of
the ${\mathcal L}$--genus
is characterized by the property
$\langle {\mathcal L}_k(N),x\rangle=\sigma(X)$ for
$x\in {\rm im}([N,S^{n-4k}] \to H_{4k}(N;\QQ))$, where
$X^{4k} \subset N$ is the inverse image $f^{-1}(p)$ of some regular value
$p\in S^{n-4k}$ of a map $f\co N\to S^{n-4k}$
which represents the Poincar\'e dual of $x$ and is $PL$ transverse regular to $p$.
Set $x'=(h^{-1})_*(x)\in H_{4k}(M;\QQ)$ and let us show that
$\langle {\mathcal L}_k(M),x'\rangle=\langle {\mathcal L}_k(N),x\rangle$.

Since $X$ is framed in $N$, it has an open $PL$ neighborhood $W=X\times \RR^{n-4k}$ in $N$.
Let $V=h^{-1}(W)\subset M$, then $h$ restricts to a homeomorphism $f\co V\to W$.
Homotope $f$ to get a map $g$ which is $PL$ transverse regular to $X=X\times\{0\}$,
and set $Y$ to be the preimage $g^{-1}(X)$,
then $\langle {\mathcal L}_k(M),x'\rangle=\sigma(Y)$ and this is equal to
$\sigma(X)=\langle {\mathcal L}_k(N),x\rangle$ by \fullref{signature}.
\end{proof}

\bibliographystyle{gtart}
\bibliography{link}

\end{document}